\documentclass[twoside,letter,final]{siamltex}
\usepackage{amsmath}
\usepackage{mathtools}
\usepackage{amssymb}
\usepackage{amsfonts}
\usepackage{amsxtra}
\usepackage{amstext}
\usepackage{amsbsy}
\usepackage{amscd}
\usepackage{graphicx}
\usepackage{float}
\usepackage{cite}
\usepackage{xcolor}
\usepackage{srcltx} 
\usepackage{marginnote}
\usepackage{rotating}
\definecolor{darkblue}{rgb}{0.0,0.0,0.6}
\definecolor{darkgreen}{rgb}{0.0,0.6,0.0}
\usepackage{pgfplots}
\pgfplotsset{width=6.5cm,compat=newest}

\usepackage{hyperref}
\usepackage[utf8]{inputenc}      
\usepackage{booktabs}            
\usepackage{url}                 
\usepackage[T1]{fontenc}         
\usepackage{algorithmic}         
\usepackage{comment}
\usepackage{calc}

\newtheorem{algorithm}{Algorithm}

\numberwithin{table}{section}    
\numberwithin{figure}{section}   
\numberwithin{equation}{section} 

\newcommand{\sign}{\operatorname{sign}}

\newcommand{\argmax}{\operatorname{argmax}}

\newcommand{\spa}{\operatorname{span}}
\newcommand{\R}{\mathbb{R}}

\newcommand{\BB}{\mathcal{B}}
\newcommand{\II}{\mathcal{I}}

\renewcommand{\AA}{\mathcal{A}}

\newcommand{\KK}{\mathcal{K}}
\newcommand{\MM}{\mathcal{M}}
\newcommand{\K}{\mathbb{K}}
\newcommand{\TT}{\mathcal{T}}
\newcommand{\Uad}{U_{\mathrm{ad}}}

\newcommand{\dual}[2]{\left\langle #1 ,\, #2 \right\rangle}

\newcommand{\ep}{\varepsilon}

\newcommand{\range}{\operatorname{range}}

\newtheorem{remark}[theorem]{Remark}

\usepackage{cleveref}

\begin{document}
\title{Bouligand Analysis and Discrete Optimal Control of Total Variation-Based Variational Inequalities}
\date{\today}
\author{J.C.~De Los Reyes\footnotemark[3]}
\renewcommand{\thefootnote}{\fnsymbol{footnote}}
\footnotetext[3]{Research Center in Mathematical Modeling and Optimization (MODEMAT), Quito, Ecuador}
\renewcommand{\thefootnote}{\arabic{footnote}}

\maketitle
\begin{abstract}
  We investigate differentiability and subdifferentiability properties of the solution mapping associated with variational inequalities (VI) of the second kind involving the discrete total-variation. Bouligand differentiability of the solution operator is established via a direct quotient analysis applied to a primal-dual reformulation of the VI. By exploiting the structure of the directional derivative and introducing a suitable subspace, we fully characterize the Bouligand subdifferential of the solution mapping. We then derive optimality conditions characterizing Bouligand-stationary and strongly-stationary points for discrete VI-constrained optimal control problems. A trust-region algorithm for solving these control problems is proposed based on the obtained characterizations, and a numerical experiment is presented to illustrate the main properties of both the solution and the proposed algorithm.
\end{abstract}


\section{Introduction}
In this paper, we continue our investigation of optimality conditions and solution algorithms for optimal control problems constrained by variational inequalities (VI) of the second kind. The inequalities considered here involve the discrete total variation (TV) seminorm and the control is of distributed nature. Such models arise in various applications, including viscoplastic fluid flow, image processing, and elastoplasticity \cite{Glowinski,scherzer2009variational,de2016optimal}.

Optimal control problems involving variational inequalities were first investigated in the late 1970s and early 1980s, with a primary focus on obstacle-type problems (see, e.g., \cite{Mignot1976,MignotPuel1984,Barbu1984,BergouniouxMignot2000}). In parallel, problems with abstract variational inequality constraints were also studied \cite{Barbu1993,BonnansTiba1991,Bergounioux1998}, leading to the derivation of general optimality conditions. However, due to the highly abstract nature of these formulations, the resulting optimality systems did not exhibit complementarity relations between the variables and lacked a precise characterization of the adjoint multipliers on the so-called biactive set.

A particular class of variational inequalities of the second kind with convex, nonsmooth, sparsity-promoting terms was studied in depth in \cite{DeLosReyes2009}, where optimality systems with complementarity relations were derived using a regularization approach. The analysis relied on a specific family of regularizing functions, yielding a limiting \emph{C-stationarity} system. More recently, a direct approach was proposed in \cite{DelosReyesMeyer2015}, focusing on the differentiability properties of the solution operator. In that work, weak directional differentiability was established for problems involving the nondifferentiable $L^1$-norm of the state, leading to the derivation of an optimality system characterizing \emph{S-stationary} points for the case of distributed controls. Extending such results to problems involving the infinite-dimensional TV-seminorm remains challenging, as it requires very restrictive assumptions on the structure of the biactive set \cite{christof2015differentiability}.

In this paper, we adopt an intermediate approach to investigate the differentiability and subdifferentiability properties of the solution operator when the nondifferentiable term in the variational inequality involves the discrete total variation seminorm. While related questions have been addressed in \cite{outrata2000,hintermullerwu2015} using tools from Mordukhovich’s generalized differentiation theory, those analyses typically rely on abstract variational principles and do not fully capture the structure of the directional derivative of the solution mapping. In contrast, our method is based on a direct quotient analysis of a primal-dual reformulation of the variational inequality. This leads to a directional differentiability result that allows us to rigorously derive both \emph{Bouligand stationarity conditions} and \emph{strong stationarity conditions} for the associated optimal control problem--going beyond the \emph{M-stationary} conditions obtained in \cite{outrata2000}.

The second goal of the paper is to analyze the subdifferential structure of the solution mapping. By introducing a suitably defined subspace, we provide a complete characterization of the Bouligand subdifferential and show that the directional derivative admits a linear representative in every direction. This result is both theoretically significant and algorithmically useful, as it underpins the design of an efficient trust-region method within the framework developed in \cite{christof2020nonsmooth}. In particular, the Bouligand subdifferential is employed to define a generalized Cauchy point based on a suitable adjoint system.

The main contributions of this work can be summarized as follows:
\begin{itemize}
    \item[i)] We develop a comprehensive analysis for variational inequalities of the second kind involving discrete total variation, providing the necessary foundation for studying nonsmooth phenomena in optimization and control.
    \item[ii)] Our approach combines a primal-dual reformulation with a direct quotient analysis to rigorously establish the Bouligand differentiability of the solution operator and to study the structure of the corresponding directional derivative.
    \item[iii)] For the first time, we provide an explicit and constructive characterization of the Bouligand subdifferential of the solution mapping associated with variational inequalities involving total variation.
    \item[iv)] The theoretical results serve as a cornerstone for deriving sharp optimality conditions for discrete optimal control problems governed by total variation-based variational inequalities, including both Bouligand and strong stationarity systems. 
    \item[v)] Our analysis further enables the design and study of nonsmooth trust-region algorithms, which critically rely on a detailed understanding of the subdifferentiability properties of the solution operator.
\end{itemize}

The structure of the paper is as follows. In Section 2, we study the directional differentiability of the solution operator associated with the variational inequality using a direct quotient analysis. We establish Bouligand differentiability, and, in the case of an empty biactive set, we also obtain Fréchet differentiability. Section 3 is devoted to characterizing the Bouligand subdifferential of the solution operator. In Section 4, we analyze the related discrete optimal control problems and derive B- and strong stationarity conditions. A trust-region algorithm is proposed in Section 5. Finally, in Section 6, we present a numerical experiment based on a Bingham flow control problem.

\section{Directional derivative of the VI solution mapping}
We are concerned with the following class of variational inequalities of the second kind: Find $y \in \R^n$ such that
\begin{equation}\label{eq:vi2}
 \dual{A y}{v-y} + \sum_{j = 1}^m \big(|(\K v)_j| - |(\K y)_j|\big) \geq \dual{u}{v-y},
 \quad \forall\, v\in \R^n,
\end{equation}
with
\begin{itemize}
 \item $A\in \R^{n\times n}$ symmetric positive definite.
 \item $K^{(i)} \in \R^{m\times n}$, $i = 1, ..., d$, discrete $i$-th partial derivative,
 \begin{equation*}
  \K : \R^n \to \R^{m\times d}, \quad \K y = (K^{(1)} y, ..., K^{(d)}y)
 \end{equation*}
 (so $\K$ is linear and bounded and thus a tensor of third order) and $\K^*: \R^{m\times d}\to \R^n$ is the adjoint mapping w.r.t.\ the scalar product associated with the Frobenius norm, i.e.,
 $\dual{A}{B}_{\R^{d\times m}}  = \sum_{i=1}^d\sum_{j=1}^m A_{ij} B_{ij}$,\\
 $(\K y)_j\in \R^d$, $j = 1, ..., m$, $j$-th row of $\K y$, corresponds to the discrete gradient at element $j$. Moreover, we assume that $\K$ is injective and, thus, the matrix $\K^* \K$ is symmetric positive definite.
 \item $|\,.\,|$ and $\dual{\,.\,}{\,.\,}$ denote the Euclidian norm and scalar product, respectively,
 in $\R^n$ as well as in $\R^d$, depending on the dimension of the corresponding input variable.
\end{itemize}

Inequality \eqref{eq:vi2} represents the necessary and sufficient optimality condition of the following strictly convex energy minimization problem
\begin{equation}\label{eq:energymin}
 \min_{y\in \R^n}\quad \frac{1}{2}\, \dual{y}{Ay} - \dual{u}{y} + \Psi(\K y).
\end{equation}
where
\begin{equation}\label{eq:defpsi}
 \psi: \R^{d} \ni w \mapsto |w| \in \R \quad \text{and} \quad
 \Psi: \R^{m\times d} \ni B \mapsto \sum_{j = 1}^m \psi(B_j) \in \R.
\end{equation}
As the objective in \eqref{eq:energymin} is uniformly convex, one readily gets the following result.

\begin{lemma} \label{lemma: global Lipschitz}
 For every $u\in \R^n$ there exists a unique solution $y\in \R^n$ of \eqref{eq:energymin} and \eqref{eq:vi2}, respectively.
 The associated solution operator $S: \R^n \ni u \mapsto y \in \R^n$ is globally Lipschitz.
\end{lemma}

By the definition of $\Psi$, \eqref{eq:vi2} is equivalent to
\begin{equation*}
 \Psi(\K v) \geq \Psi(\K y) + \dual{u - Ay}{v - y}
 \quad \Longleftrightarrow\quad u - Ay \in \partial (\Psi \circ \K)(y) = \K^* \partial \Psi(\K y),
\end{equation*}
where we used the chain rule for convex subdifferentials since $\Psi$ is convex and continuous.
Thus, there exists a dual multiplier $q\in \partial \Psi(\K y)$ such that $u - Ay = \K^* q$, which results in
\begin{subequations}\label{eq:complsys}
\begin{align}
 & A y + \K^* q = u\label{eq:complsysa}\\
 & \dual{q_j}{(\K y)_j} = |(\K y)_j|, && \forall \, j = 1, ..., m,\label{eq:complsysb}\\ 
 & |q_j| \leq 1 && \forall \, j = 1, ..., m,\label{eq:complsysc}
\end{align}
\end{subequations}
where $q_j \in \R^d$, $j = 1, ..., m$, denotes $j$-th row of $q$. Let us define the active and inactive sets by
\begin{equation}\label{eq:activeset}
 \II(y) := \{ j\in \{1, ..., m\} : (\K y)_j \neq 0\} \quad\text{and}\quad \AA(y) := \{1, ..., m\} \setminus \II(y),
\end{equation}
and the biactive set by 
\begin{equation}\label{eq:biactiveset}
 \BB(y) := \{ j\in \{1, ..., m\} : |q_j| = 1 \wedge (\K y)_j = 0\}.
\end{equation}
Then \eqref{eq:complsysb} yields that $q$ satisfies
\begin{equation}\label{eq:quni}
 q_j = \frac{(\K y)_j}{|(\K y)_j|}, \quad \forall\, j \in \II(y),
\end{equation}
so that the components of $q$ in $\II(y)$ are uniquely determined by $y$.
Note that, in general, $q$ need not be unique on the set $\AA(y)$.

In all what follows, we call a vector $q$ satisfying \eqref{eq:complsys} \emph{slack variable}. Moreover, the argument in the active, inactive and biactive sets notation will be omitted if it can be clearly inferred from the context.

\begin{lemma}\label{lem:ky}
 Let $y\in \R^n$ and $q\in \R^{m\times d}$ be given. Then the set $\KK(y)$ defined by
 \begin{equation}\label{eq:Ky}
 \begin{aligned}
  \KK(y) := \{v\in \R^n: \; & (\K v)_j = 0, \text{ if } |q_{j}| < 1,\\
  &\dual{q_{j}}{(\K v)_j} = |(\K v)_j|, \text{ if } |q_{j}| = 1 \wedge (\K y)_j = 0 \}
 \end{aligned}
 \end{equation}
 is a convex cone. If $y$ and $q$ satisfy \eqref{eq:complsysb}-\eqref{eq:complsysc}, this set can equivalently be expressed as
 \begin{equation}\label{eq:coneequiv}
 \begin{aligned}
  \KK(y) = \Bigg\{v\in \R^n :
  \dual{\K^*q}{v} \geq \sum_{j\in \II(y)} \Big\langle \frac{(\K y)_j}{|(\K y)_j|} , (\K v)_j \Big\rangle
  + \sum_{j \in \AA(y)} |(\K v)_j| \Bigg\}.
 \end{aligned}
 \end{equation}
\end{lemma}

\begin{proof}
 Thanks to $|q_{j}| = 1$ and the Cauchy-Schwarz inequality, the last condition \eqref{eq:Ky} is
 equivalent to
 \begin{equation*}
  \dual{q_{j}}{(\K v)_j} \geq |(\K v)_j|, \quad \text{if }\; |q_{j}| = 1 \wedge (\K y)_j = 0.
 \end{equation*}
 Then the linearity of $\K$ and the convexity of $|\,.\,|$ immediately yield the first result.

 To proof the equivalent reformulation in case that $q$ and $y$ satisfy \eqref{eq:complsysb}-\eqref{eq:complsysc}, denote the set in \eqref{eq:coneequiv} by $\MM$.
 Thanks to \eqref{eq:quni} and the definition of $\KK(y)$ in \eqref{eq:Ky} we immediately obtain $\KK(y)\subset \MM$.
 To proof the converse inclusion, let $v\in \MM$ be arbitrary. Then \eqref{eq:quni} implies
 \begin{equation*}
 \begin{aligned}
 \dual{q}{\K v}_{\R^{m\times d}}
  =  \sum_{j\in \II(y)} \Big\langle \frac{(\K y)_j}{|(\K y)_j|} , (\K v)_j \Big\rangle
  + \sum_{j \in \AA(y)} \dual{q_j}{(\K v)_j},
 \end{aligned}
 \end{equation*}
 and, consequently, by the definition of $\MM$,
 \begin{equation*}
  \sum_{j \in \AA(y)} |(\K v)_j| \leq \sum_{j \in \AA(y)} \dual{q_j}{(\K v)_j} \leq \sum_{j \in \AA(y)} |(\K v)_j|,
 \end{equation*}
 where we used the Cauchy-Schwarz inequality and $|q_j| \leq 1$, see \eqref{eq:complsysc}, for the last estimate.
 Thus we obtain
 \begin{equation*}
  0 = \sum_{j \in \AA(y)} \Big(|(\K v)_j| - \dual{q_j}{(\K v)_j} \Big).
 \end{equation*}
 Again due to the Cauchy-Schwarz inequality and $|q_j| \leq 1$, every addend in the above sum is non-negative so that
 \begin{equation}\label{eq:activeeq}
  |(\K v)_j| = \dual{q_j}{(\K v)_j}, \quad \forall\, j \in \AA(y),
 \end{equation}
 is obtained. Since by \eqref{eq:complsysb} there holds
 \begin{equation*}
 \begin{aligned}
  \AA(y) &= \{ j \in \{1, ..., m\} : (\K y)_j = 0\}\\
  &= \{ j \in \{1, ..., m\} : |q_j| < 1\} \cup \{ j \in \{1, ..., m\} : |q_j| = 1 \wedge (\K y)_j = 0\},
 \end{aligned}
 \end{equation*}
 \eqref{eq:activeeq} finally yields that $v\in \KK(y)$.
\end{proof}

\begin{remark}\label{rem:coneequal}
 The above lemma shows the following:
 If $q^1$ and $q^2$ are two different slack variables associated with the solution $y$ of \eqref{eq:vi2},
 then the two sets
 \begin{equation*}
 \begin{aligned}
  \KK_i :=  \{v\in \R^n: \; & (\K v)_j = 0, \text{ if } |q_{j}^i| < 1,\\
  &\dual{q_{j}^i}{(\K v)_j} = |(\K v)_j|, \text{ if } |q_{j}^i| = 1 \wedge (\K y)_j = 0 \}, \quad i = 1,2,
 \end{aligned}
 \end{equation*}
 coincide, since $\K^* q^1 = u - A y = \K^* q^2$. Therefore, the set in \eqref{eq:coneequiv} is the same in both cases.
 This also justifies the notation $\KK(y)$, as this set does not depend on the slack variable, but only on the solution $y$.
\end{remark}

Next, let $h\in \R^n$ be given and consider the perturbed problem
\begin{equation}\label{eq:vi2pert}
 \dual{A y^t}{v-y^t} + \sum_{j = 1}^m \big(|(\K v)_j| - |(\K y^t)_j|\big) \geq \dual{u + t \,h}{v-y^t},
 \quad \forall\, v\in \R^n.
\end{equation}
The Lipschitz continuity of $S$ readily yields
\begin{equation*}
 \Big|\frac{y^t - y}{t}\Big| \leq c\,|h|,
\end{equation*}
and, hence, a subsequence of $\{(y^t - y)/t\}$ converges to some $\eta \in \R^n$. Without loss of generality, we denote this subsequence by the same symbol, i.e.,
\begin{equation}\label{eq:convy}
 \frac{y^t - y}{t} \to \eta.
\end{equation}
As before one can reformulate the VI in terms of a complementarity system, i.e.,
\begin{subequations}\label{eq:pert}
\begin{gather}
 A y^t + \K^* q^t = u + t\,h \label{eq:perta}\\
 \dual{q_{j}^t}{(\K y^t)_j} = |(\K y^t)_j|,\quad |q_{j}^t| \leq 1 \quad \forall \, j = 1, ..., m.\label{eq:pertb}
\end{gather}
\end{subequations}
In view of \eqref{eq:pertb}, the sequence $\{q^t\}$ is bounded and therefore a subsequence, again w.l.o.g.\ denoted by the same symbol,
exists so that
\begin{equation}\label{eq:convq}
 q^t \to \tilde q \in \R^n.
\end{equation}
Due to \eqref{eq:convy}, we additionally have $y^t \to y$ such that we can pass to the limit $t\searrow 0$ in \eqref{eq:pert} to obtain
\begin{gather*}
 A y + \K^* \tilde q = u\\
 \dual{\tilde q_j}{(\K y)_j} = |(\K y)_j|,\quad |\tilde q_j| \leq 1 \quad \forall \, j = 1, ..., m,
\end{gather*}
such that $\tilde q$ belongs to the set of slack variables associated with $y$.
This in particular implies that \eqref{eq:quni} holds with $q = \tilde q$.

\begin{proposition}\label{lem:etaink}
  It holds that $\eta \in \KK(y)$.
\end{proposition}

\begin{proof}
 Adding the complementarity relations in \eqref{eq:complsysb} and \eqref{eq:pertb} gives
 \begin{equation} \label{eq: increments in q inactive}
  \Big\langle q^t_j \, ,\, \frac{(\K y^t)_j - (\K y)_j}{t} \Big\rangle
  + \frac{1}{t}\, \dual{q^t_j - q_j}{(\K y)_j} = \frac{|(\K y^t)_j| - |(\K y)_j|}{t}.
 \end{equation}
 Now let $j\in \AA(y)$ be arbitary so that $(\K y)_j = 0$. In this case the above equation becomes
 \begin{equation}\label{eq:compldiff}
  \Big\langle q^t_j \, ,\, \frac{(\K y^t)_j - (\K y)_j}{t} \Big\rangle
  = \frac{|(\K y^t)_j| - |(\K y)_j|}{t}
 \end{equation}
 and, thanks to \eqref{eq:convy}, \eqref{eq:convq}, and the Bouligand differentiability of
 $\psi: \R^d \ni v \mapsto |v| \in \R$,
 we can pass to the limit in \eqref{eq:compldiff} to obtain
 \begin{equation*}
  \dual{\tilde q_j}{(\K\eta)_j} = \psi'\big((\K y)_j;(\K\eta)_j\big) = |(\K\eta)_j|, \quad \forall\, j\in \AA(y).
 \end{equation*}
 Arguing as at the end of the proof of Lemma \ref{lem:ky}, cf.\ \eqref{eq:activeeq},
 and keeping Remark \ref{rem:coneequal} in mind (note that $\tilde q$ is a slack variable), we get that $\eta \in \KK(y)$.
\end{proof}

%

\begin{lemma}\label{lem:sign1}
 For every $v\in \KK(y)$ there holds
 \begin{equation*}
  \Big\langle \K^*\, \frac{q^t - q}{t}\,,\, v \Big\rangle \leq
  \sum_{j\in \II(y)}  \frac{1}{t}\, \Big\langle\frac{(\K y^t)_j }{|(\K y^t)_j|} - \frac{(\K y)_j }{|(\K y)_j|}
  \,,\,(\K v)_j \Big\rangle,
 \end{equation*}
 for all $t>0$ sufficiently small.
\end{lemma}

\begin{proof}
 Let $v\in \KK(y)$ be arbitrary.
 Due to $y^t \to y$, there holds $\II(y) \subset \II(y^t)$, provided that $t>0$ is sufficiently small.
 Hence, $q_j^t = (\K y^t)_j / |(\K y^t)_j|$ in $\II(y)$, cf.\ \eqref{eq:quni}, giving in turn
 \begin{equation*}
 \begin{aligned}
  \dual{\K^* q^t}{v}
  & = \sum_{j\in \II(y)} \dual{q^t_j}{(\K v)_j} + \sum_{j \in \AA(y)} \dual{q^t_j}{(\K v)_j}\\
  & \leq \sum_{j\in \II(y)} \Big\langle \frac{(\K y^t)_j}{|(\K y^t)_j|} , (\K v)_j \Big\rangle
  + \sum_{j \in \AA(y)} |q^t_j| |(\K v)_j|.
 \end{aligned}
 \end{equation*}
 Employing $|q_j^t|\leq 1$, see \eqref{eq:pertb}, and the second formulation of $\KK(y)$ in \eqref{eq:coneequiv}
 implies the result.
\end{proof}

\begin{lemma}\label{lem:sign2}
 For all $t > 0$ and all $j\in \{1, ..., m\}$, there holds
 \begin{equation*}
  \Big\langle \frac{q^t_j - q_j}{t}\,,\, \frac{(\K y^t)_j - (\K y)_j}{t} \Big\rangle \geq 0.
 \end{equation*}
\end{lemma}
\begin{proof}
 The complementarity relations in \eqref{eq:complsysb} and \eqref{eq:pertb} yield
 \begin{equation*}
 \begin{aligned}
  &\Big\langle \frac{q^t_j - q_j}{t}\,,\, \frac{(\K y^t)_j - (\K y)_j}{t} \Big\rangle \\
  &\quad = \frac{1}{t^2} \Big( \dual{q^t_j}{(\K y^t)_j}
  - \dual{q^t_j}{(\K y)_j} - \dual{q_j}{(\K y^t)_j} + \dual{q_j}{(\K y)_j} \Big)\\
  &\quad \geq \frac{1}{t^2} \Big( |(\K y^t)_j| - \underbrace{|q^t_j|}_{\leq 1} |(\K y^t)_j|
  - \underbrace{|q_j|}_{\leq 1}|(\K y)_j| + |(\K y)_j| \Big) \geq 0.
 \end{aligned}
 \end{equation*}
\end{proof}

\begin{theorem}
 The solution operator $S: \R^n \to \R^n$ associated with \eqref{eq:vi2} is directionally differentiable. Its
 directional derivative at $u\in \R^n$, in direction $h\in \R^n$, is the unique solution $\eta\in \R^n$ of the following
 VI of the first kind
 \begin{equation}\label{eq:vi1}
 \left.
 \begin{aligned}
  & \eta \in \KK(y),\\
  & \dual{A\eta}{v-\eta} +
  \begin{aligned}[t]
   \sum_{j\in \II(y)} \Big\langle\frac{(\K \eta)_j }{|(\K y)_j|}
   - \dual{(\K y)_j}{(\K\eta)_j} \,\frac{(\K y)_j}{|(\K y)_j|^3} \,,\,(\K v)_j - (\K\eta)_j \Big\rangle
  \end{aligned}\\[1mm]
  &\qquad \qquad \quad \geq \dual{h}{v - \eta} \quad \forall\, v\in \KK(y),
 \end{aligned}
 \right\}
 \end{equation}
 where $y = S(u)$, $\II(y)$ and $\AA(y)$ are the sets defined in \eqref{eq:activeset} and $\KK(y)$ is given by
 \begin{equation*}
 \begin{aligned}
  \KK(y) = \Bigg\{v\in \R^n :
  \dual{u - A y}{v} \geq \sum_{j\in \II(y)} \Big\langle \frac{(\K y)_j}{|(\K y)_j|} , (\K v)_j \Big\rangle
  + \sum_{j \in \AA(y)} |(\K v)_j| \Bigg\}.
 \end{aligned}
 \end{equation*}
\end{theorem}

\begin{proof}
 First, the condition $\eta \in \KK(y)$ was already proven in Lemma \ref{lem:etaink}.
 To verify the VI in \eqref{eq:vi1}, let $v\in \KK(y)$ be arbitrary.
 We test \eqref{eq:complsysa} and \eqref{eq:perta} with $v - (y^t - y)/t$ and subtract the arising equations.
 In this way we obtain, for all $t>0$ sufficiently small, the following estimate, by using Lemmata \ref{lem:sign1} and \ref{lem:sign2},
 \begin{equation}\label{eq:vi1est}
 \begin{aligned}
  & \Big\langle h \,,\, v - \frac{y^t - y}{t} \Big\rangle -
  \Big\langle A \,\frac{y^t - y}{t} \,,\, v - \frac{y^t - y}{t} \Big\rangle \\
  &\quad = \Big\langle \K^*\, \frac{q^t - q}{t}\,,\, v \Big\rangle
  \begin{aligned}[t]
   & - \sum_{i\in \II(y)} \Big\langle \frac{q^t_j - q_j}{t}\,,\, \frac{(\K y^t)_j - (\K y)_j}{t} \Big\rangle\\
   & - \sum_{i\in \AA(y)} \Big\langle \frac{q^t_j - q_j}{t}\,,\, \frac{(\K y^t)_j - (\K y)_j}{t} \Big\rangle
  \end{aligned}\\
  &\quad \leq \sum_{j\in \II(y)}  \Big\langle \frac{1}{t}\Big( \frac{(\K y^t)_j }{|(\K y^t)_j|} - \frac{(\K y)_j }{|(\K y)_j|} \Big)
  \,,\,(\K v)_j - \frac{(\K y^t)_j - (\K y)_j}{t} \Big\rangle.
 \end{aligned}
 \end{equation}
 Note that, for $t>0$ sufficiently small, we have $\II(y) \subset \II(y^t)$ and thus $q_j^t = (\K y^t)_j / |(\K y^t)_j|$ in $\II(y)$,
 which was already used in the proof of Lemma \ref{lem:sign1}. As $\psi$, defined in \eqref{eq:defpsi},
 is smooth on $\R^d\setminus\{0\}$, its derivative given by
 \begin{equation*}
  \nabla \psi(w) = \frac{w}{|w|},\quad  w\in \R^d\setminus\{0\}
 \end{equation*}
 is differentiable at $(\K y)_j$, for all $j\in \II(y)$.
 Together with \eqref{eq:convy}, this allows to pass to the limit in \eqref{eq:vi1est}, which, in view of
 \begin{equation*}
  \psi''(w) = \frac{1}{|w|}\, I - \frac{1}{|w|^3}\, w\,w^\top  \quad \forall \,w\in \R^d\setminus \{0\},
 \end{equation*}
 implies \eqref{eq:vi1}. Thus we have shown that the limit $\eta$ of a subsequence of $\{(y^t - y)/t\}_{t>0}$ satisfies \eqref{eq:vi1}.

 To prove the convergence of the whole sequence, we just have to show that the limit is unique.
 For this purpose, observe that \eqref{eq:vi1} is the necessary optimality condition of the following minimization problem
 \begin{equation}\label{eq:derivmin}
  \left.
  \begin{aligned}
   \min_{\eta \in \R^d} ~ & f_y(\eta) := \frac{1}{2}\, \dual{\eta}{A\eta} - \dual{h}{\eta}
   + \frac{1}{2} \sum_{i\in \II(y)} \Big(\frac{|(\K\eta)_j|^2}{|(\K y)_j|} - \frac{\dual{(\K y)_j}{(\K\eta)_j}^2}{|(\K y)_j|^3}\Big)\\
   \text{s.t.}\quad & \eta \in \KK(y).
  \end{aligned}
  \right\}
 \end{equation}
 The feasible set $\KK(y)$ is convex by Lemma \ref{eq:coneequiv}. For the second derivative of the objective, the Cauchy-Schwarz inequality and the coercivity of $A$ yield
 \begin{equation*}
  w^\top f_y''(\eta) w = w^\top A w + \sum_{i\in \II(y)} \frac{1}{|(\K y)_j|} \Big( |(\K w)_j|^2
  - \frac{\dual{(\K y)_j}{(\K w)_j}^2}{|(\K y)_j|^2}\Big) > 0,
 \end{equation*}
 for all $w\in \R^n\setminus\{0\}$, so that the objective in \eqref{eq:derivmin} is strictly convex.
 Thus \eqref{eq:derivmin} is a stricly convex minimization problem and consequently \eqref{eq:vi1} is
 also sufficient for optimality and thus equivalent to \eqref{eq:derivmin}.
 The strict convexity yields the uniqueness of the solution $\eta$, which finally finishes the proof.
\end{proof}

\begin{remark}
  Using the definitions of $\psi$ and $\Psi$ in \eqref{eq:defpsi}, the VI in \eqref{eq:vi1} can equivalently be written in short form as
 \begin{equation*}
 \left.
 \begin{aligned}
  & \eta \in \KK(y),\\
  & \dual{A\eta}{v-\eta} +
  \sum_{j\in \II(y)} (\K\eta)_j^\top \psi''((\K y)_j) (\K(v-\eta))_j
  \geq \dual{h}{v - \eta} \quad \forall\, v\in \KK(y),
 \end{aligned}
 \quad\right\}
 \end{equation*}
 with $\KK(y) = \{v\in \R^n : \dual{u - A y}{v} \geq \Psi'(\K y;\K v) \}.$
\end{remark}

\begin{remark}\label{rem:bouli}
 As $S$ is globally Lipschitz continuous, its directional differentiability automatically implies that
 $S$ is \underline{Bouligand-differentiable} (see, e.g., \cite[Thm.\ 3.1.2]{Scholtes-Buch}).
\end{remark}

\begin{corollary} \label{cor: Gateaux dif}
 If there exists a slack variable $q$ such that the strict complementarity condition
 \begin{equation}\label{eq:strict}
  (\K y)_j = 0 \quad \Longrightarrow \quad |q_j| < 1
 \end{equation}
 holds true, then the directional derivative $\eta$ solves the following \emph{linear system}:
 \begin{equation}\label{eq:Sabl}
 \begin{aligned}
  & A \eta + \K^* \lambda = h,\\
  & \lambda_j - \frac{(\K \eta)_j }{|(\K y)_j|} + \dual{(\K y)_j}{(\K\eta)_j} \,\frac{(\K y)_j}{|(\K y)_j|^3} = 0, & &\forall\, j\in\II(y),\\
  & (\K \eta)_j = 0, & & \forall \, j\in\AA(y),
 \end{aligned}
 \end{equation}
 with a slack variable $\lambda \in \R^{m\times d}$. The solution operator $S:\R^n \to \R^n$ of \eqref{eq:vi2} is therefore
 Fr\'echet differentiable in case that \eqref{eq:strict} holds.
\end{corollary}

\begin{proof}
 If there is a slack variable such that \eqref{eq:strict} holds, then, according to Lemma \ref{lem:ky}, the convex cone $\KK(y)$ becomes
 \begin{equation}\label{eq:subspace}
  \KK(y) = \{v\in \R^n : (\K v)_j = 0, \text{whenever } (\K y)_j = 0\},
 \end{equation}
 and, consequently, $\KK(y)$ is a linear subspace in this case. The VI in \eqref{eq:vi1} thus becomes a variational equation
 so that the directional derivative of $S$ is a linear mapping w.r.t.\ the direction $h$.
 Since $S$ is Bouligand-differentiable, see Remark \ref{rem:bouli}, this yields the Fr\'echet-differentiability.
 To derive the precise form of the derivative in \eqref{eq:Sabl}, consider again the minimization problem
 \eqref{eq:derivmin}, which is equivalent to the VI in \eqref{eq:vi1}.
 If $\KK(y)$ takes the form \eqref{eq:subspace}, then the KKT-conditions for this problem look as follows:
 \begin{equation}\label{eq:kkt}
 \left.
 \begin{gathered}
  \begin{aligned}[t]
   \dual{A \eta - h}{v}
   + \sum_{j\in \II(y)} \Big\langle\frac{(\K \eta)_j }{|(\K y)_j|}
   - \dual{(\K y)_j}{(\K\eta)_j} \,\frac{(\K y)_j}{|(\K y)_j|^3} \,,\,(\K v)_j \Big\rangle \qquad &\\[-1mm]
   + \sum_{j\in \AA(y)} \dual{\nu_j}{(\K v)_j} = 0 &
  \end{aligned}\\
  (\K \eta)_j = 0 \quad \forall\, j \in \AA(y).
 \end{gathered}
 \quad \right\}
 \end{equation}
 with Lagrange-multipliers $\nu_j\in \R^d$, $j\in \AA(y)$. Note that the Abadie constraint qualification is satisfied, since
 the constraints of \eqref{eq:derivmin} are linear such that \eqref{eq:kkt} is necessary and, due to convexity, sufficient
 for optimality. If we introduce $\lambda \in \R^{m\times d}$ by
 \begin{equation*}
  \lambda_j :=
  \begin{cases}
   \nu_j, & j \in \AA(y),\\
   \frac{(\K \eta)_j }{|(\K y)_j|} - \dual{(\K y)_j}{(\K\eta)_j} \,\frac{(\K y)_j}{|(\K y)_j|^3}, & j \in \II(y),
  \end{cases}
 \end{equation*}
 then \eqref{eq:Sabl} is obtained.
\end{proof}

The above corollary suggests the following algorithm to verify strict complementarity and compute the Fr\'echet derivative of $S$:
\begin{algorithm}
\begin{algorithmic}[1]
 \STATE Solve \eqref{eq:vi2} to obtain $y$
 \STATE Compute a slack variable $\hat q$ (if it is not a by-product of Step~1.)
 \IF{\eqref{eq:strict} is fulfilled with $q = \hat q$}
  \STATE Compute the derivative $\eta = S'(u)h$ by solving the linear system \eqref{eq:Sabl}.
 \ELSE
 \STATE Solve the following minimization problem (with $y$ from Step~1):
  \begin{equation}\label{eq:minq}
   \begin{aligned}
    \min_{q\in \R^{m\times d}, \, r\in \R} &\quad \frac{1}{2}\, r^2\\
    \text{subject to: } & \quad
    \begin{aligned}[t]
     & A y + \K^* q = u,\\
     & \dual{q_j}{(\K y)_j} = |(\K y)_j|, \\ 
     & |q_j|^2 \leq 1 \quad \forall\, j = 1, ..., m,\\
     & |q_j|^2 \leq r \quad \forall\, j\in \AA(y).
    \end{aligned}
   \end{aligned}
  \end{equation}
  with solution $\bar r$ and $\bar q$
  \IF{$\bar r < 1$}
   \STATE Compute the derivative $\eta = S'(u)h$ by solving the linear system \eqref{eq:Sabl}.
  \ELSE
   \STATE $S$ is not Fr\'echet-differentiable at $u$.
  \ENDIF
 \ENDIF
\end{algorithmic}
\end{algorithm}

By solving the optimization problem \eqref{eq:minq}, one computes the slack variable with the minimum $\ell_\infty$-norm on $\AA(y)$.
Thus, if there is a slack variable satisfying \eqref{eq:strict}, it will be detected by solving \eqref{eq:minq}.

\section{Bouligand subdifferential}
We now focus on the study of the Bouligand subdifferential of the solution operator $S(u)$ and obtain a linear system of equations that characterizes its elements.

\begin{theorem} \label{thm: characterization of B-subdif}
   Let $G$ be an element of $\partial_B S(u)$ and let $y=S(u)$ be the solution of \eqref{eq:vi2}. There exists a partition $\BB_0 \cup \BB_1$ of the biactive set $\BB$ such that, for any $h \in \R^n$, $G h=:\tilde \eta \in V$ corresponds to the unique solution of the system
\begin{subequations} \label{eq: Bouligand subdiff. system}
  \begin{align}
    & \dual{A \tilde \eta}{v} + \sum_{j \in \mathcal I} \left \langle \tilde \lambda_j , (\K v)_j \right \rangle = \dual{h}{v}, &&  \text{for all } v \in V\\
     & \tilde \lambda_j = \frac{(\K \tilde \eta)_j }{|(\K y)_j|} - \dual{(\K y)_j}{(\K \tilde \eta)_j} \,\frac{(\K y)_j}{|(\K y)_j|^3}  && \text{for }j\in\II.
  \end{align}
\end{subequations}
where $V:=\{ v \in \R^n: (\K v)_j=0, \forall j \in \AA_s \cup \BB_0; (\K v)_j \in \spa (q_j), \forall j \in \BB_1 \}$ and $\AA_s:= \{ j: |q_j| <1 \}$.
\end{theorem}
\begin{proof}
  Let $D_S \subset \mathbb R^n$ denote the set where $S$ is differentiable. By definition of the Bouligand subdifferential, there is a sequence $\{u_n\} \subset D_S$ such that $u_n \to u$ and $S'(u_n) \to G$. Thanks to the Lipschitz continuity of $S$, we know that
  \begin{equation*}
    y_n=S(u_n) \to S(u):=y \quad \text{ and } \quad \K^* q_n = u_n - A y_n \to u-Ay = \K^* q.
  \end{equation*}
  The last representation follows from the fact that $\{q_n\}$ is also bounded and has therefore a convergent subsequence. The claim follows from the uniqueness of the limit.

  Considering the inactive and strongly active sets:
  \begin{equation*}
    \II= \{ j: (\K y)_j \not = 0 \}, \qquad \AA_S=\{ j: |q_j| <1 \},
  \end{equation*}
  it follows by continuity that $\II \subset \II^n$ and $\AA_S \subset \AA_S^n$, for $n \geq N$ sufficiently large, where $\II^n$ and $\AA_S^n$ correspond to the inactive and strongly active sets associated to $u_n$. Since $\{u_n\} \subset D_S$, it then follows, for $h \in \R^n$, that $\eta_n := S'(u_n)h$ satisfies the system (see \eqref{eq:Sabl})
  \begin{align}
    &A \eta_n + \K^* \lambda_n = h, &&\\
    &(\lambda_n)_j - \frac{(\K \eta_n)_j }{|(\K y_n)_j|} + \dual{(\K y_n)_j}{(\K \eta_n)_j} \,\frac{(\K y_n)_j}{|(\K y_n)_j|^3} = 0 && j \in \II^n,\\
    & (\K \eta_n)_j=0, && j \in \AA^n,
  \end{align}
  or, equivalently,
  \begin{subequations} \label{eq: Gateaux diff. system}
      \begin{multline}
        \dual{A \eta_n}{v} + \sum_{j \in \II^n} \left \langle \frac{(\K \eta_n)_j }{|(\K y_n)_j|} - \dual{(\K y_n)_j}{(\K \eta_n)_j} \,\frac{(\K y_n)_j}{|(\K y_n)_j|^3},(\K v)_j \right \rangle\\
        = \dual{h}{v}, \,\text{ for all } v \in V_n,
      \end{multline}
      \begin{equation}
        (\K \eta_n)_j=0, \hspace{2cm} j \in \AA_n,
      \end{equation}
  \end{subequations}
where $V_n:=\{ v \in \R^n: (\K v)_j=0, \forall j \in \AA_n \}$.
From the definition of the Bouligand subdifferential it follows that $\tilde \eta = \lim_{n \to +\infty} \eta_n.$ Moreover, since for $j \in \II$ the sequence $\{ (\lambda_n)_j \}$ is bounded, there is a convergent subsequence with a limit $\tilde \lambda_j$. Consequently, up to a subsequence, by passing to the limit we get that
  \begin{align}
    &A \tilde \eta + \K^* \tilde \lambda = h &&\\
    &\tilde \lambda_j - \frac{(\K \tilde \eta)_j }{|(\K y)_j|} + \dual{(\K y)_j}{(\K \tilde \eta)_j} \,\frac{(\K y)_j}{|(\K y)_j|^3} = 0 && j\in\II\\
    & (\K \tilde \eta)_j=0, && j \in \AA_s.
  \end{align}

It remains to analyze what happens on the biactive set $\BB= \{ j: (\K y)_j = 0, |q_j|=1 \}.$ Let us first consider the subset
$$\BB_0=\{ j \in \BB: \exists \text{ a subsequence }\{y_{n_k}\}: (\K y_{n_k})_j=0, \forall k \}.$$ Since $\eta_n \to \tilde \eta,$ we get that
\begin{equation*}
  (\K \tilde \eta)_j =0, \quad \text{ for all } j \in \AA_S \cup \BB_0.
\end{equation*}

Considering now the subset
\begin{equation*}
  \BB_1:= \BB \backslash \BB_0 = \{ j \in \BB: (\K y_n)_j \not =0, \forall n \in \mathbb N \text{ suff. large} \},
\end{equation*}
and since $j \in \II^n$, we obtain for any $v \in V$ that
\begin{align*}
  \dual{\tilde \lambda_j}{(\K v)_j}
  &= \lim_{n \to +\infty} \dual{(\lambda_n)_j}{(\K v)_j}= \lim_{n \to +\infty} (c_n)_j \dual{(\lambda_n)_j}{(q_n)_j}\\
  &=  \lim_{n \to +\infty} (c_n)_j \dual{(\lambda_n)_j}{\frac{(\K y_n)_j}{|(\K y_n)_j|}}\\
  &=  \lim_{n \to +\infty} \frac{(c_n)_j}{|(\K y_n)_j|^2} \dual{\left(I - \frac{(\K y_n)_j (\K y_n)_j^T}{|(\K y_n)_j|^2} \right) (\K \tilde \eta_n)_j}{(\K y_n)_j}\\
  &=  \lim_{n \to +\infty} \frac{(c_n)_j}{|(\K y_n)_j|^2} \dual{\left(I - \frac{(\K y_n)_j (\K y_n)_j^T}{|(\K y_n)_j|^2} \right) (\K y_n)_j}{(\K \tilde \eta_n)_j}\\
  &=0.
\end{align*}
Passing to the limit in equation \eqref{eq: Gateaux diff. system} then yields \eqref{eq: Bouligand subdiff. system}.

Finally, we prove that, for $j \in \BB_1$, $(\K \tilde \eta)_j \in \spa(q_j)$. To do so, note that, thanks to \eqref{eq: Gateaux diff. system} and the positive definitness of $A$, we obtain, testing the equation with $v = \eta_n$, that
  \begin{align*}
    0 & \leq \frac{|(\K \eta_n)_j|^2}{|(\K y_n)_j|} - \frac{\dual{(\K y_n)_j}{(\K \eta_n)_j}^2}{|(\K y_n)_j|^3} \\ 
    & \leq \dual{A \eta_n}{\eta_n} + \sum_{j \in \II^n} \left \langle \frac{(\K \eta_n)_j }{|(\K y_n)_j|} - \dual{(\K y_n)_j}{(\K \eta_n)_j} \,\frac{(\K y_n)_j}{|(\K y_n)_j|^3},(\K \eta_n)_j \right \rangle\\
    & = \dual{h}{\eta_n},
  \end{align*}
Since $\{ \eta_n \}$ is bounded, there exists a constant $C>0$ such that
  \begin{equation*}
    0 \leq \frac{1}{|(\K y_n)_j|} \left( |(\K \eta_n)_j|^2- \frac{\dual{(\K y_n)_j}{(\K \eta_n)_j}^2}{|(\K y_n)_j|^2} \right) \leq C, \quad \text{for } j \in \BB_1.
  \end{equation*}
Since $(\K y_n)_j \to 0$, we conclude that
  \begin{equation*}
    0 = \lim_{n \to \infty} |(\K \eta_n)_j|^2 - \frac{\dual{(\K y_n)_j}{(\K \eta_n)_j}^2}{|(\K y_n)_j|^2} = |(\K \tilde \eta)_j|^2 - \dual{q_j}{(\K \tilde \eta)_j}^2,
  \end{equation*}
which implies, since $|q_j|=1$, that $(\K \tilde \eta)_j = c_j q_j$ for some $c_j \in \mathbb R$. Consequently, $(\K \tilde \eta)_j \in \spa(q_j)$ and the proof is complete.
\end{proof}

\begin{corollary}\label{cor: G is a solution of the linear system}
  Let $G\in\partial_B S(u)$. There exists a partition of the biactive set $\BB = \BB_0 \cup \BB_1$ and a multiplier $\theta \in \R^{m \times d}$ such that, for any $h$, $\tilde{\eta}:= Gh$ is the unique solution of the system
  \begin{subequations}\label{eq:full G_linear_system}
    \begin{align}
      & A \tilde{\eta} + \K^T \theta = h\\
      & (\K \tilde{\eta})_j  =0, && \forall j \in \AA_s \cup \BB_0 \label{eq:full G_linear_system_2}\\ 
      & (\K \tilde{\eta})_j  \in \spa (q_j), && \forall j \in \BB_1 \label{eq:full G_linear_system_3}\\
      & \theta_j  = \frac{(\K \tilde \eta)_j }{|(\K y)_j|} - \dual{(\K y)_j}{(\K \tilde \eta)_j} \,\frac{(\K y)_j}{|(\K y)_j|^3}, &&\forall j\in\II,\\
      &\dual{\theta_j}{q_j} =0, &&\forall j\in\BB_1.
    \end{align}
  \end{subequations}
\end{corollary}
\begin{proof}
  Consider the functional defined by
  \begin{equation*}
    \dual{\mathcal F}{v}:= \dual{A \tilde{\eta}}{v} + \sum_{j\in\II} \dual{\frac{(\K \tilde \eta)_j }{|(\K y)_j|} - \dual{(\K y)_j}{(\K \tilde \eta)_j} \,\frac{(\K y)_j}{|(\K y)_j|^3}}{(\K v)_j} - \dual{h}{v},
  \end{equation*}
for all $v\in V$. It is clear that system \eqref{eq: Bouligand subdiff. system} can equivalently be written as $\mathcal F \in V^\perp$. Moreover, the linear subspace $V$, can be represented as
$$
    V= \left( \bigcap_{j \in \AA_S \cup \BB_0} V^0_j \right) \cap \left( \bigcap_{j \in \BB_1} V^1_j \right),
$$
  where
  \begin{align*}
    & V^0_j:= \{ v \in \R^n: (\K v)_j =0 \}, &&  j \in \AA_S \cup \BB_0,\\
    & V^1_j:= \{ v \in \R^n: (\K v)_j \in \spa(q_j) \}, &&  j \in \BB_1.
  \end{align*}
  It then follows that the orthogonal complement of $V$ can be expressed as $V^\perp = \sum_{j \in \AA_S \cup \BB_0} (V^0_j)^\perp + \sum_{j \in \BB_1} (V^1_j)^\perp.$

  For $j \in \AA_S \cup \BB_0$, we readily obtain that $(V^0_j)^\perp= \ker (\K_j)^\perp$ and, thanks to the orthogonality relations, also $\ker (\K_j)^\perp=\range (\K_j^\top)$. Consequently, for any $\xi_j \in (V^0_j)^\perp$, there is a $\pi_j \in \R^2$ such that $\xi_j= \K_j^\top \pi_j$, and
  \begin{equation*}
    \sum_{j \in \AA_S \cup \BB_0} (V^0_j)^\perp = \sum_{j \in \AA_S \cup \BB_0} \K_j^\top \pi_j.
  \end{equation*}

  Any element $v \in V^1_j$, with $j \in \BB_1$, can be represented as sum of an element from the nullspace and the row space of $\K_j$, i.e.,
  \begin{equation*}
    v= \phi + \varphi, \quad \text{ with }(\K_j \varphi)=0 \text{ and } \phi \in \range (\K_j^\top).
  \end{equation*}
  Since $(\K v)_j \in \spa(q_j)$ and $(\K_j \varphi)=0$, it also follows that $(\K \phi)_j \in \spa(q_j)$. Taking an element $w_j \in (V^1_j)^\perp$, it can be represented as $w_j=\tilde{w}_j+\hat{w}_j$, with $\tilde{w}_j \in \range(\K_j^\top)$ and $\hat{w}_j \in \range(\K_j^\top)^\perp = \ker(\K_j)$. Consequently, there exists $\psi_j$ such that
  \begin{equation*}
    w_j = \K_j^\top \psi_j + \hat{w}_j, \quad \text{ with }\K_j \hat{w}_j=0.
  \end{equation*}
  Multiplying $w_j$ with $v_j \in  V^1_j$ we get, for some $\sigma \in \mathbb R^{m\times d}$,
  \begin{align*}
    (w_j,v_j)  &= (\K_j^\top \psi_j + \hat{w}_j, \phi + \varphi)\\ 
    &= \dual{\psi_j}{\K_j \phi} +  (\hat{w}_j, \K_j^\top \sigma) + (\hat{w}_j, \varphi)\\
    & = c \dual{\psi_j}{q_j} + (\hat{w}_j, \varphi),
  \end{align*}
  since $\K_j \varphi =\K_j \hat{w}_j= 0$. For the product to be zero, it is then required that $(\hat{w}_j, \varphi)=0, \forall \varphi \in \ker(\K_j)$, and $\dual{\psi_j}{q_j}=0.$
  Since $\hat{w}_j$ belongs to $\ker(\K_j)$ as well, it follows that $\hat{w}_j=0.$
  Thus,
  \begin{equation*}
    \sum_{j \in \BB_1} (V^1_j)^\perp = \sum_{j \in \BB_1} \K_j^\top \psi_j, \quad \psi_j \in \R^2: \dual{\psi_j}{q_j}=0.
  \end{equation*}

  Altogether, we get existence of multipliers $\pi_j$ and $\psi_j$ such that
  \begin{equation*}
    \mathcal F + \sum_{j \in \AA_S \cup \BB_0} \K_j^\top \pi_j + \sum_{j \in \BB_1} \K_j^\top \psi_j =0,
  \end{equation*}
  with $\dual{\psi_j}{q_j}=0$. Defining
\begin{equation*}
  \theta_j := \begin{cases}
    \frac{(\K \tilde \eta)_j }{|(\K y)_j|} - \dual{(\K y)_j}{(\K \tilde \eta)_j} \,\frac{(\K y)_j}{|(\K y)_j|^3} , & j\in\II,\\
    \pi_j, & j \in \AA_S \cup \BB_0,\\
    \psi_j , & j \in \BB_1,
\end{cases}
\end{equation*}
we obtain the desired result. The system \eqref{eq:full G_linear_system} is then equivalent to \eqref{eq: Bouligand subdiff. system} and, since $G$ is an element of the Bouligand subdifferential, it follows that $\tilde \eta = Gh$ is a solution of the linear system \eqref{eq:full G_linear_system}.
\end{proof}


We consider next the converse implication and prove that for any splitting of the biactive set $\BB$, the corresponding solution $\tilde \eta = G h$ of system \eqref{eq: Bouligand subdiff. system} characterizes an element of the Bouligand subdifferential $\partial_B S(u)$.

\begin{theorem} \label{thm: evry linear soliution is element of bouligand}
Let $\tilde \eta = G h$ be a solution of system \eqref{eq: Bouligand subdiff. system} for a given partition $\BB=\BB_0 \cup \BB_1$. Then $G$ corresponds to an element of $\partial_B S(u)$.
%
\end{theorem}
\begin{proof}
Let $\BB_0 \subset \BB$ be arbitrary but fix and $\BB_1 = \BB \backslash \BB_0$. Without loss of generality, we assume that $(\K \tilde \eta)_j = c_j q_j \not = 0$, for all $j \in \BB_1$. Otherwise we may consider the modified set $\widetilde \BB_0 = \BB_0 \cup \{ j \in \BB_1: (\K \tilde \eta)_j=0\}$
and the corresponding equivalent system \eqref{eq: Bouligand subdiff. system}.

We will next show that there exists a sequence $\{ u_n \}$ such that
\begin{align*}
  &u_n \in D_S, \qquad (\K y_n)_j=0, ~\forall j \in \AA_S \cup \BB_0, \qquad (\K y_n)_j \not =0, ~\forall j \in \II \cup \BB_1,\\
  &\text{and} \quad u_n \to u, \quad S'(u_n) \to G, \quad \text{as }n \to \infty.
\end{align*}

Let $\{\ep_n\} \subset \mathbb R_+$ be a sequence such that $\ep_n \to 0$ as $n\to \infty$, and consider a sequence $\{y_{n}\}$ such that 
$$(\K y_n)= (\K y) + \ep_n C (\K \tilde \eta),$$
where $C$ is a diagonal matrix with 
$$
  C_{jj} = 
    \begin{cases}
      c_j^{-1} & \text{if }j \in \BB_1,\\
      1 & \text{otherwise},
    \end{cases}
$$ 
where $c_j \in \mathbb R$ is the constant arising from \eqref{eq:full G_linear_system_2}. Existence of such sequence can be obtain thanks to the invertibility of $\K^* \K.$

For $j \in \II$ it then follows that
$$|(\K y_n)_j|= |(\K y)_j + \ep_n (\K \tilde \eta)_j|,$$
which implies that $|(\K y_n)_j| \not =0,$ for $\ep_n>0$ sufficiently small. On the other hand, $|(\K y_n)_j| = \ep_n |c_j^{-1}| |(\K \tilde \eta)_j| \not = 0$, for $j \in \BB_1$. Consequently, $\II^{n}= \II \cup \BB_1$ and, thanks to \eqref{eq:full G_linear_system_2}, $\AA^{n}= \AA \backslash \BB_1$.

On $\II^n$ we define the multiplier $q_j^n= \frac{(\K y^n)_j}{|(\K y^n)_j|}$, for $j \in \II^n$, which implies that $q_j^n= \frac{c_j^{-1} (\K \tilde \eta)_j}{|c_j^{-1}||(\K \tilde \eta)_j|}= \frac{q_j}{|q_j|}= q_j$, for $j \in \BB_1.$

On the set $\BB_0$ we define
$$i_j^* = \argmax_{i \in \{1, \dots, d\}} |q_{ji}|$$
and consider the canonical vectors
\begin{equation*}
  (e_j^*)_i= \begin{cases} 0 &\text{ if }i \not = i_j^*,\\ 1 &\text{ if }i = i_j^*,
\end{cases} \qquad j \in \BB_0
\end{equation*}
Moreover, we consider the perturbed multiplier
$$q_j^n = q_j - \ep_n \sign(q_{j i_j^*}) e^*_j.$$
It then follows that
$$|q_j^n| = |q_j - \ep_n \sign(q_{j i_j^*}) e^*_j| <|q_j|=1, ~j \in \BB_0.$$
Taking $q_j^n:=q_j$, for $j \in \AA_S$, we then get that $\AA_S^n =\AA_S \cup \BB_0 = \AA \backslash \BB_1$, which implies that $\BB^n=\AA^n \backslash \AA^n_S= \emptyset.$ Moreover, it can be verified that $|q_j^n| \leq 1, \forall j,$ and, for $j \in \BB_1$, we get that
$$\dual{q_j^n}{(\K y^n)_j}=\dual{\frac{(\K \tilde \eta)_j}{|(\K \tilde \eta)_j|}}{\ep_n (\K \tilde \eta)_j}= \ep_n |(\K \tilde \eta)_j|=|(\K y^n)_j|.$$
The sequence $\{ q^n \}$ converges therefore to the dual multiplier $q$, since $q_j^{n} \to q_j$, for $j \in \II \cup \BB_0$, and $q_j^n = q_j$, for $j \in \AA_S \cup \BB_1$.

Introducing $\xi=\frac{1}{\ep_n}(q^n-q)$ and using the control
$$u^n=u+ \ep_n A \tilde \eta+ \ep_n \K^* \xi$$
it then follows that
\begin{align*}
  & A y^n + \K^* q^n = u^n,\\
  & \dual{q^n_j}{(\K y^n)_j} = |(\K y^n)_j|, \forall j\\
  & |q^n_j| \leq 1 \quad \forall j.
\end{align*}
Since $\BB^n = \emptyset$, we get that $u^n \in D_S$ and, moreover, $u^n \to u$ as $n \to \infty$.

It remains to verify that $S'(u^n) \to G$. Thanks to the Lipschitz continuity of $S$ we get that, for $\ep_n \to 0$,
$$ \| S'(u^{n}) \| \leq L, \qquad \forall n.$$
Therefore, there exists a subsequence $\{u^{n_k}\}$ and a limit $H \in \mathbb R^{n \times n}$ such that $S'(u^{n_k}) \to H \in \partial_B S(u)$, as $k \to \infty.$ Since system \eqref{eq: Bouligand subdiff. system} is uniquely solvable, the result $H=G$ follows from the uniqueness of the limit.
\end{proof}

As a consequence of the previous two results, we may obtain a characterization of the generalized jacobian of the solution mapping as well. This is the content of the following corollary.
\begin{corollary}
  An element $G$ belongs to the generalized jacobian $\partial S(u)$ if and only if, for any $h \in \R^n$, $G h=:\hat \eta \in \hat V$ corresponds to the unique solution of the system
\begin{subequations} \label{eq: Clarke subdiff. system}
 \begin{align}
   & \dual{A \hat \eta}{v} + \sum_{j \in \mathcal I} \left \langle \hat \lambda_j , (\K v)_j \right \rangle = \dual{h}{v}, &&  \text{for all } v \in \hat V\\
    & \hat \lambda_j = \frac{(\K \hat \eta)_j }{|(\K y)_j|} - \dual{(\K y)_j}{(\K \hat \eta)_j} \,\frac{(\K y)_j}{|(\K y)_j|^3}  && \text{for }j\in\II.
 \end{align}
\end{subequations}
where $\hat V:=\{ v \in \R^n: (\K v)_j=0, \forall j \in \AA_s; (\K v)_j \in \spa (q_j), \forall j \in \BB \}.$
\end{corollary}

Next we verify that, along a given direction, there exists a solution of the linear system \eqref{eq: Bouligand subdiff. system}, which coincides with the directional derivative. When properly characterized, this enables the use of a linear representative of the (otherwise nonlinear) directional derivative within any solution algorithm (see Section \ref{sec: TR algorithm} below).
\begin{theorem}
For any $u,h \in \R^n$, there exists a linearized element $\tilde \eta =G h$, solution of \eqref{eq: Bouligand subdiff. system}, such that $S'(u;h)=G h$.
\end{theorem}
\begin{proof}
Let us recall that the directional derivative of the solution operator, in direction $h$, is given by the unique solution $\eta \in \KK(y)$ of
\begin{equation}\label{eq: directional derivative in proof of Thm 5.1}
\left.
\begin{aligned}
 & \dual{A\eta}{v-\eta} +
 \begin{aligned}[t]
  \sum_{j\in \II(y)} \Big\langle\frac{(\K \eta)_j }{|(\K y)_j|}
  - \dual{(\K y)_j}{(\K\eta)_j} \,\frac{(\K y)_j}{|(\K y)_j|^3} \,,\,(\K v)_j - (\K\eta)_j \Big\rangle
 \end{aligned}\\[1mm]
 &\qquad \qquad \quad \geq \dual{h}{v - \eta}, \quad \forall\, v\in \KK(y),
\end{aligned}
\right\}
\end{equation}
where $\KK(y)$ is given by \eqref{eq:Ky}. Defining the matrices $T_j:=\frac{1}{|(\K y)_j|} \left( I- \frac{(\K y)_j (\K y)_j^T}{|(\K y)_j|^2}\right),$ for $j \in \II(y)$, and the linear operator $L: \R^n \to \R^n$ such that, for $w \in \R^n$,
\begin{equation*}
      \dual{Lw}{v}:=\dual{Aw}{v}+ \sum_{j \in \II(y)}\dual{T_j(\K w)_j}{(\K v)_j}, \quad \forall v \in \R^n,
\end{equation*}
inequality \eqref{eq: directional derivative in proof of Thm 5.1} can be expressed as
\begin{equation*}
      \dual{L \eta}{v- \eta} \geq \dual{h}{v-\eta}, \quad \forall v \in \KK(y)
\end{equation*}
or, equivalently, as $\eta = P_{\KK}(\eta-\sigma (L\eta+ h))$, for all $\sigma>0$, where $P_{\KK}$ stands for the projection onto the convex cone $\KK(y)$.

Let us now consider the sets $\BB_0:=\{ j \in \BB: (\K \eta)_j=0\}$ and $\BB_1= \BB \backslash \BB_0$. Since $\eta \in \KK(y)$, it follows that $(\K \eta)_j= c_j q_j$, for all $j \in \BB_1$, for some $c_j > 0$. Therefore, $\eta$ belongs to the subspace
$$V:=\{ v \in \R^n: (\K v)_j=0, \forall j \in \AA_s \cup \BB_0; (\K v)_j \in \spa (q_j), \forall j \in \BB_1 \}.$$


Additionally, for any $w \in V$ it follows that $\eta \pm t w \in \KK(y)$, for $t$ sufficiently small. Using these vectors in \eqref{eq: directional derivative in proof of Thm 5.1} then yields
\begin{equation*}
  \dual{A\eta}{w}+ \sum_{j \in \II(y)}\dual{T_j(\K \eta)_j}{(\K v)_j}= \dual{h}{w}, \quad \forall w \in V.
\end{equation*}
Therefore, the directional derivative takes the form $\eta = Gh$, solution of \eqref{eq: Bouligand subdiff. system}, with $\BB_0$ and $\BB_1$ as defined above.
\end{proof}


%


\section{Stationarity conditions}
We focus next on the study of optimality conditions for the discrete (VI)-constrained optimal control problem:
\begin{subequations} \label{VI constrained optimization problem}
\begin{align}
\min_{u \in U_{ad}} ~& J(y,u)\\
\text{subject to: }& \langle Ay,v-y \rangle +|\K v|_1- |\K y|_1 \geq \langle u, v-y\rangle,  \text{ for all } v \in \mathbb R^n,
\end{align}
\end{subequations}
where we assume that $J$ is continuously differentiable, $U_{ad}$ is a closed convex set, and $A$ and $\K$ are defined as in equation \eqref{eq:vi2}. The goal along this section will be the characterization of stationary points for problem \eqref{VI constrained optimization problem}, through a system of necessary optimality conditions that include properties of the adjoint state on the biactive set.

By using the solution operator $S(u)$ of the variational inequality, the problem can be reformulated in reduced form as
\begin{equation} \label{eq:reduced optimization problem}
  \min_{u \in U_{ad}} ~f(u)=J(S(u),u).
\end{equation}
Thanks to the chain rule for B-differentiable functions (see, e.g., \cite[Section~4.1]{cui2021modern}), it follows that the composite mapping $f$, as a function of $u$, is B-differentiable as well. The directional derivative is given by $$f'(u;h)= \nabla_y J(S(u),u)^T \eta+ \nabla_u J(S(u),u)^T h,$$ with $\eta \in \R^n$ the unique solution to \eqref{eq:vi1}. Moreover, if $\bar u$ is a local optimal solution, then it satisfies the following necessary condition:
\begin{equation} \label{eq: necessary condition}
  f(\bar u;u-\bar u)=\nabla_y J(\bar y,\bar u)^T \bar \eta+ \nabla_u J(\bar y,\bar u)^T (u- \bar u) \geq 0, \text{ for all } u \in U_{ad},
\end{equation}
where $\bar y := S(\bar u)$ and $\bar \eta$ corresponds to the solution to \eqref{eq:vi1} with $h=u- \bar u$. A point $\bar u$ satisfying the necessary condition \eqref{eq: necessary condition} is called \underline{B-stationary}.

Let us next consider, for a given $u \in U_{ad}$, the tangent cone
\begin{align*}
  \TT(u) := \left\{ (\eta, h) : \exists \{u_n\} \subset U_{ad}, \{t_n\} \subset \R^+ \text{ s.t. } \frac{u_n-u}{t_n} \to h, \frac{S(u_n)-S(u)}{t_n} \to \eta \right\}.
\end{align*}

\begin{theorem} \label{thm: B-stationarity}
  Let $\bar u \in U_{ad}$ be a local optimal solution of \eqref{VI constrained optimization problem} and $\bar y = S(\bar u)$. Then $\bar u$ satisfies the following inequality:
  \begin{equation}
    \nabla_y J(\bar y,\bar u)^T \eta + \nabla_u J(\bar y,\bar u)^T h \geq 0, \text{ for all } (\eta,h) \in \TT(\bar u).
  \end{equation}
\end{theorem}

\begin{proof}
  Let $(\eta,h) \in \TT(\bar u)$. From the definition of the tangent cone, there exist sequences $\{u_n\} \subset U_{ad}$ and $\{t_n\} \subset \R^+$ such that $\frac{u_n-u}{t_n} \to h$ and $\frac{S(u_n)-S(u)}{t_n} \to \eta$. From \eqref{eq: necessary condition} and the positive homogeneity of the Bouligand derivative it follows that
  \begin{equation}
    \nabla_y J(\bar y,\bar u)^T S' \left( \bar u; \frac{u_n- \bar u}{t_n} \right) + \nabla_u J(\bar y,\bar u)^T \left( \frac{u_n- \bar u}{t_n} \right) \geq 0.
  \end{equation}
  Thanks to the Lipschitz continuity of the B-derivative of $S$ with respect to the direction and the continuous differentiability of $J$, we may pass to the limit in the previous inequality and get the result.
\end{proof}

For the case $U_{ad}= \R^n$ we are able to obtain a multiplier characterization of local minima, which leads to a \underline{strong stationarity} optimality system.

\begin{theorem}
  Let $\bar u$ be a local optimal solution of \eqref{VI constrained optimization problem}, with $U_{ad}= \R^n$, and $\bar y = S(\bar u)$. Then there exist multipliers $p \in \R^n$ and $\mu \in \R^n$ such that
\begin{subequations}
  \begin{align}
     & A y + \K^* q = u\\
     & \dual{q_j}{(\K y)_j} = |(\K y)_j|, && \forall j = 1, ..., m\\ 
     & |q_j| \leq 1, && \forall j = 1, ..., m\\
     & \dual{A p}{v}+ \sum_{j\in \II(\bar y)} \dual{T_j (\K p)_j}{(\K v)_j} = \dual{\nabla_y J(\bar y, \bar u)- \mu}{v}, && \forall v \in \R^n\\
     & p \in \KK(\bar y)\\
     & \dual{\mu}{\phi}  \geq 0, && \forall \phi \in \KK(\bar y)\\
     & p + \nabla_u J(\bar y, \bar u)=0,
  \end{align}
\end{subequations}
where $T_j:= \frac{1}{|(\K \bar y)_j|} \left( I- \frac{(\K \bar y)_j (\K \bar y)_j^T}{|(\K \bar y)_j|^2}\right),$ for $j \in \II(\bar y).$
\end{theorem}
\begin{proof}
  Let us define the projection operator $P : \R^n \to \KK(\bar y)$ which assigns to each $\xi \in \R^n$ the unique $P(\xi)$ solution of
  \begin{equation*}
    a(P(\xi), \phi- P(\xi)) \geq a(\xi, \phi-P(\xi)), \quad \forall \phi \in \KK(\bar y),
  \end{equation*}
  where $a(\cdot, \cdot)$ is the coercive bilinear form defined by
  $$ a(v,w):= \dual{Av}{w}+ \sum_{j \in \II(\bar y)}\dual{T_j (\K v)_j}{(\K w)_j}, \quad \forall v, w \in \R^n.$$
  Moreover, we denote by $L$ the symmetric positive matrix associated with $a(\cdot, \cdot)$, i.e., $\dual{L v}{w}:= a(v,w), \, \forall v,w \in \R^n$.

  The polar cone of $\KK(\bar y)$ with respect to $a(\cdot, \cdot)$ is given by
  $$\left( \KK(\bar y) \right)^0_a:= \{ \varphi \in \R^n: a(\varphi, \phi) \leq 0, \quad \forall \phi \in \KK(\bar y) \}.$$ By defining $Q(\xi)= \xi- P(\xi)$, it can be easily verified that $Q(\xi) \in \left( \KK(\bar y) \right)^0_a$ and, moreover, $a(Q(\xi),P(\xi))=0.$

  With help of these operators, the Bouligand derivative of the solution mapping can be written as
  $$S'(\bar u; h)= P(L^{-1} h),$$ 
  since $\dual{h}{\phi}= a(L^{-1}h, \phi),$ for all $\phi \in \KK(\bar y).$ Consequently, the directional derivative of the cost function can be written as
  \begin{align*}
    f'(\bar u; h)&= \dual{\nabla_y J(\bar y, \bar u)}{S'(\bar u;h)}+\dual{\nabla_u J(\bar y, \bar u)}{h}\\
    & = a(L^{-1} \nabla_y J(\bar y, \bar u), P(L^{-1}h))+ a(L^{-1}h, \nabla_u J(\bar y, \bar u))\\
    & = a(P(L^{-1}h), L^{-1} \nabla_y J(\bar y, \bar u))+ a(P(L^{-1}h)+Q(L^{-1}h), \nabla_u J(\bar y, \bar u))\\
    &= a(P(L^{-1}h), L^{-1} \nabla_y J(\bar y, \bar u)+ \nabla_u J(\bar y, \bar u))+ a(Q(L^{-1}h), \nabla_u J(\bar y, \bar u)).
  \end{align*}

Defining $\xi_0:= -L^{-1} \nabla_y J(\bar y, \bar u)-\nabla_u J(\bar y, \bar u)$ and $\xi_1:= - \nabla_u J(\bar y, \bar u)$ we then get that
  \begin{align*}
    f'(\bar u; h)&=-a(P(L^{-1}h), \xi_0) -a(Q(L^{-1}h), \xi_1)\\
    &=-a(P(L^{-1}h), P(\xi_0))- a(P(L^{-1}h), Q(\xi_0))\\
    & \hspace{1cm} -a(Q(L^{-1}h), P(\xi_1)) -a(Q(L^{-1}h), Q(\xi_1)).
  \end{align*}
For the choice $h_0= L P(\xi_0) \Leftrightarrow L^{-1} h_0 = P(\xi_0)$ we obtain that
\begin{align*}
  & P(L^{-1}h_0)= P(P(\xi_0))=P(\xi_0)\\
  & Q(L^{-1}h_0)= Q(P(\xi_0))=0.
\end{align*}
Consequently, from the B-stationarity condition \eqref{eq: necessary condition}, we get that
\begin{equation*}
  f'(\bar u; h_0) = - a(P(\xi_0),P(\xi_0))-a(Q(\xi_0),Q(\xi_0))=-a(P(\xi_0),P(\xi_0)) \geq 0,
\end{equation*}
which implies that $P(\xi_0)=0$ or, equivalently, $\xi_0 \in \left( \KK(\bar y) \right)^0_a$.

On the other hand, for the choice $h_1= L Q(\xi_1) \Leftrightarrow L^{-1} h_1 = Q(\xi_1)$ we obtain that
\begin{align*}
  & P(L^{-1}h_1)= P(Q(\xi_1))=0\\
  & Q(L^{-1}h_1)= Q(Q(\xi_1))=Q(\xi_1).
\end{align*}
Therefore,
\begin{equation*}
  f'(\bar u; h_1) = - a(Q(\xi_1),P(\xi_1))-a(Q(\xi_1),Q(\xi_1))=-a(Q(\xi_1),Q(\xi_1)) \geq 0,
\end{equation*}
and, thus, $Q(\xi_1)=0$ or, equivalently, $\xi_1=P(\xi_1) \in \KK(\bar y)$.

Defining $\mu:= - L \xi_0$ and the adjoint state $p:= L^{-1}(\nabla_y J(\bar y, \bar u)- \mu)=\xi_1$, we then get that $p+ \nabla_u J(\bar y, \bar u)=0$ and
\begin{equation*}
  - a(\xi_0, \phi)= -a(L^{-1} \mu, \phi)=\dual{\mu}{\phi} \geq 0, \quad \forall \phi \in \KK(\bar y),
\end{equation*}
which concludes the proof.
\end{proof}


\section{A nonsmooth trust region algorithm} \label{sec: TR algorithm}
In this section we divise a trust-region algorithm for solving \eqref{eq:reduced optimization problem}. Due to the nonsmoothness of the problem, we consider a quadratic model involving an element of the Bouligand subdifferential, instead of the cost function gradient. However, this choice alone may not lead to a convergent sequence of iterates (see, e.g., \cite{apkarian2016nonsmooth}), as Cauchy points do not take neighborhood information into account. To ensure convergence, we introduce an additional phase in the algorithm, triggered when the trust-region radius becomes small, in which a generalized model is considered (see \cite{christof2020nonsmooth} for further details).



Let us start by describing the first phase of the algorithm. As shown previously (see Corollary \ref{cor: Gateaux dif}), in the case of an empty biactive set, additional differentiability properties of the solution mapping may be obtained. Indeed, in this case, the derivative is of Fréchet type and is characterized by \eqref{eq:Sabl}. Based on this expression, the existence of a classical adjoint state can be established, allowing for the application of adjoint calculus. 

Whenever the biactive set is not empty, however, the characterization of the Bouligand subdifferential enables us to introduce a \emph{generalized adjoint} state associated to system \eqref{eq:full G_linear_system}. To do so, let us consider a partition $\BB_0 \cup \BB_1$ of the biactive set and define the adjoint state $p \in \mathbb R^n$ as the solution to the system:
\begin{align}\label{eq: adjoint system}
& Ap + \K^* \lambda = \nabla_y J(y,u), &&\\
& \lambda_j = \frac{(\K p)_j}{|(\K y)_j|}- \frac{(\K y)_j (\K y)_j^T }{|(\K y)_j|^3}(\K p)_j, && \forall j \in \II,\\
& (\K p)_j=0, && \forall j \in \AA_S \cup \BB_0,\\
&(\K p)_j \in \spa(q_j), &&\forall j\in\BB_1,\\
&\dual{\lambda_j}{q_j} =0, &&\forall j\in\BB_1.
\end{align}
With this generalized adjoint at hand, we may consider the corresponding Bouligand subdifferential of the cost function as follows:
\begin{equation} \label{eq: generalized gradient}
\partial_B f(u) \ni g= \nabla_y J(y,u) + p.
\end{equation}
Other elements of $\partial_B f(u)$ corresponding to different splittings of the biactive set $\BB$ may be considered as well.

Let us remark that the slack multiplier $q \in \mathbb R^{m \times d}$ is not necessarily unique, which may lead to different biactive sets and, therefore, different (and possibly unstable) numerical behavior. To remedy this, we consider hereafter the choice of the slack multiplier with the smallest Euclidean norm.

Using \eqref{eq: generalized gradient}, a quadratic model of the reduced cost function is then given by
\begin{equation}\label{eq:trmodel}
 \mathfrak{q}_k(s)=f(u_k)+ g_k^T s+ \frac{1}{2} s^T H_k s,
\end{equation}
where $g_k \in \partial_B f(u)$ and $H_k$ is a matrix with curvature information, obtained for instance with some variant of the BFGS method. The trust region radius is denoted by $\Delta_k$ and the actual and predicted reductions are defined by
$$\operatorname{ared}_k(s^k):=f(u_k)-f(u_k+s^k) \quad \text{ and } \quad \operatorname{pred}_k(s^k)=f(u_k)-\mathfrak{q}_k(s^k),$$
respectively. The quality indicator in the first phase is computed by $$\rho_k(s^k)=\frac{\operatorname{ared}_k(s^k)}{\operatorname{pred}_k(s^k)}.$$

For the second phase of the algorithm, when $\Delta_k$ is smaller than a threshold radius $\Delta_{min}$, we first identify the set of possible bi-active indices
\begin{align*}
	\mathcal P(u_k, \Delta_k) &:= \{ i \in \{1, \dots, m\}: |(\K y(u_k))_i| \leq L_y \Delta_k \land |q_i(u_k)| \geq 1- L_y \Delta \},\\
  \mathcal A_v (u_k, \Delta_k) &:= \{ i \in \{1, \dots, m\}: |q_i(u_k)| < 1- L_y \Delta \},
\end{align*}
where $L_y$ stands for the Lipschitz constant of the solution mapping. Denoting the subsets of $\mathcal P(u_k, \Delta_k)$ by $\BB^k_1, ..., \BB^k_{m_k}$ with $m_k = 2^{|\mathcal P(u_k, \Delta_k)|}$, we consider the quadratic model
\begin{equation}\label{eq:trmodelgeneralized}
 \mathfrak{q}_k(s)=f(u_k)+ \zeta+ \frac{1}{2} s^T H_k s,
\end{equation}
where $\zeta$ has to satisfy the inequalities
\begin{equation*}
  \dual{g^k_j}{d} \leq \zeta, \qquad\forall\,  j = 1, ..., m_k.
\end{equation*}
An alternative quality indicator has to be considered in this case, which is given by
\begin{equation*}
 \rho_k :=
 \begin{cases}
  \displaystyle{\frac{f(u_k) - f(u_k + d_k)}{f(u_k) - {q}_k(d_k)}},
  & \text{if }  \psi(u_k, \Delta_k) > \|g_k\|\, \Delta_k \\
  0, & \text{if }  \psi(u_k, \Delta_k) \leq \|g_k\|\, \Delta_k.
 \end{cases}
\end{equation*}

The resulting trust region algorithm is given through the following steps:\\

\begin{algorithm}[Trust-Region Algorithm for the solution of \eqref{VI constrained optimization problem}]\label{alg:trvi}
 \begin{algorithmic}[1]
  \STATE Initialization: Choose constants
  \begin{equation*}
   \Delta_{\min} > 0, \quad 0 < \eta_1 < \eta_2 < 1, \quad 0 < \beta_1 < 1 < \beta_2, \quad 0 < \mu \leq 1
  \end{equation*}
  an initial value $u_0 \in \R^n$, and an initial TR-radius $\Delta_0 > \Delta_{\min}$. Set $k=0$.
  \REPEAT
   \STATE Choose a subset $\BB_k \subseteq \BB(u_k)$, solve the generalized \emph{adjoint equation}
   \begin{align*}
     & A p_k + \K^* \lambda_k = \nabla_y J(y_k, u_k), &&\\
     & (\lambda_k)_j = \frac{(\K p_k)_j}{|(\K y_k)_j|}- \frac{(\K y_k)_j (\K y_k)_j^T }{|(\K y_k)_j|^3}(\K p_k)_j, && ~j \in \II(u_k),\\
     & (\K p_k)_j=0, && ~j \in \AA_S(u_k) \cup \BB_k,
   \end{align*}
   and set $g_k = p_k + \nabla_u J(y_k, u_k)$.
   \STATE Choose a matrix $H_k \in \R^{n\times n}_{\textup{sym}}$.
   \IF{$g_k = 0$}
    \STATE STOP the iteration, $0 \in \partial_B f(u_k)$.
   \ELSE
    \IF{$\Delta_k > \Delta_{\min}$}
     \STATE Compute an inexact solution $d_k$ of the \emph{trust-region subproblem}
     \begin{equation}\tag{Q$_{k}$}
      \left.
      \begin{aligned}
       \min_{d\in \R^n} & \quad \mathfrak{q}_k(d) := f(u_k) + \dual{g_k}{d} + \frac{1}{2}\, d^\top H_k d\\
       \text{s.t.} & \quad |d| \leq \Delta_k,
      \end{aligned}
      \quad\right\}
     \end{equation}
     that satisfies the \emph{generalized Cauchy-decrease condition}
     \begin{equation*}
      f(u_k) - \mathfrak{q}_k(d_k) \geq
      \frac{\mu}{2}\,|g_k|\,\min\Big\{ \Delta_k, \frac{|g_k|}{|H_k|} \Big\}.
     \end{equation*}
     \STATE Compute the quality indicator
     \begin{equation*}
      \rho_k := \frac{f(u_k) - f(u_k + d_k)}{f(u_k) - \mathfrak{q}_k(d_k)}.
     \end{equation*}
    \ELSIF{$\Delta_k \leq \Delta_{\min}$}
     \STATE Identify the set of possibly bi-active indices $\mathcal P(u_k, \Delta_k)$
     and their subsets $\BB^k_1, ..., \BB^k_{m_k}$.
     \FOR{$i = 1, ..., m_k$}
      \STATE Solve the \emph{adjoint equation}
      \begin{align*}
        & A p_i^k + \K^* \lambda_i^k = \nabla_y J(y_k, u_k), &&\\
        & (\lambda_i^k)_j = \frac{(\K p_i^k)_j}{|(\K y_k)_j|}- \frac{(\K y_k)_j (\K y_k)_j^T }{|(\K y_k)_j|^3}(\K p_i^k)_j, && ~j \in \II(u_k),\\
        & (\K p_i^k)_j=0, && ~j \in \AA_S(u_k) \cup \BB_i^k
      \end{align*}
      and set $g^k_j = p^k_j + \nabla_u J(y_k, u_k)$.
     \ENDFOR
      \STATE Compute an inexact, but feasible solution $d_k$ of the \emph{modified trust-region subproblem}
      \begin{equation}\tag{$\mathfrak{Q}_k$}\label{eq:qkvimod}
       \left.
       \begin{aligned}
        \min_{\zeta \in \R, d\in \R^n}  \quad & \mathfrak{q}_k(d,\zeta)
        := f(u_k) + \zeta + \frac{1}{2}\, d^\top H_k d\\
        \text{\textup{s.t.}}  \quad & |d| \leq \Delta_k, \\
        & \dual{g^k_j}{d} \leq \zeta \quad\forall\,  j = 1, ..., m_k .
       \end{aligned}
       \quad\right\}
      \end{equation}
      that satisfies the \emph{modified Cauchy-decrease condition}
     \begin{equation}\label{eq:cauchyvi}
      f(u_k) - \mathfrak{q}_k(d_k, \zeta_k) \geq
      \frac{\mu}{2}\,\psi(u_k, \Delta_k) \,\min\Big\{ \Delta_k, \frac{\psi(u_k, \Delta_k)}{\|H_k\|} \Big\}.
     \end{equation}
		 where $\psi= - \min_{|d| \leq 1} \left\{ \xi: \dual{g^k_j}{d} \leq \xi, ~\forall j=1, \dots, m_k\right\}.$
     \STATE Compute the modified quality indicator
     \begin{equation*}
      \rho_k :=
      \begin{cases}
       \displaystyle{\frac{f(u_k) - f(u_k + d_k)}{f(u_k) - {q}_k(d_k)}},
       & \text{if }  \psi(u_k, \Delta_k) > \|g_k\|\, \Delta_k \\
       0, & \text{if }  \psi(u_k, \Delta_k) \leq \|g_k\|\, \Delta_k.
      \end{cases}
     \end{equation*}
    \ENDIF
    \STATE Update: Set
    \begin{align*}
     u_{k+1} & :=
     \begin{cases}
      u_k, & \text{if } \rho_k \leq \eta_1 \quad \text{(null step)},\\
      u_k + d_k, & \text{otherwise} \quad \text{(successful step)},
     \end{cases} \\
     \Delta_{k+1} &:=
     \begin{cases}
      \beta_1\,\Delta_k, & \text{if } \rho_k \leq \eta_1,\\
      \max\{\Delta_{\min}, \Delta_k\}, & \text{if } \eta_1 < \rho_k \leq \eta_2,\\
      \max\{\Delta_{\min}, \beta_2 \Delta_k\}, & \text{if } \rho_k > \eta_2.
     \end{cases}
    \end{align*}
    Set $k = k+1$.
   \ENDIF
  \UNTIL{$0\in \partial f(u_k)$}.
 \end{algorithmic}
\end{algorithm}


For the computation of the inexact step in the previous algorithm (step 9.), we consider a dogleg strategy, which is described next. The main purpose of this choice is to accelerate the behaviour of the trust-region method, although no theoretical guarantee is available.
\begin{algorithm}{(Choice of Cauchy point)}
\begin{algorithmic}[1]
\STATE Choose the parameter values $0<\eta_1<\eta_2<1$, $0<\gamma_0<\gamma_1<1<\gamma_2$, $\Delta_{min}\geq0.$
\STATE Compute the Cauchy step $s_c^k=-t^* g_k,$ where
$$t^*=\left\{
\begin{matrix}
\displaystyle\frac{\Delta_k}{|g_k|},\ \ \text{ if }\ g_k^\top H_kg_k\leq0\\
\\
\min\left(\displaystyle\frac{|g_k|^2}{g_k^\top H_kg_k},\displaystyle\frac{\Delta_k}{|g_k|}\right),\ \ \text{ if }\ g_k^\top H_k g_k>0
\end{matrix}
\right.$$
and the Newton step $s^k_n= -H_k^{-1} g_k$.
\IF{$s^k_n$ satisfies the fraction of Cauchy decrease:
$$\exists \delta \in \,]0,1] \text{ and } \beta \geq 1 \text{ such that }|s^k| \leq \beta \Delta_k \text{ and }pred_k(s^k) \geq \delta ~pred_k(s_c^k).$$} \vspace{0.1cm}
\STATE $s^k=s^k_n$,\vspace{0.1cm}
\ELSE \vspace{0.1cm}
\STATE $s^k=s^k_c$.\vspace{0.1cm}
\ENDIF
\end{algorithmic}
\end{algorithm}

\section{Numerical experiment}
In this section we experimentally verify some properties of the proposed trust-region algorithm by means of the discretized viscoplastic Bingham flow control problem \cite{DeLosReyes2009,dlRe2010}. We focus particularly on:
\begin{itemize}
\item Total iteration number with respect to the Tikhonov regularization parameter;
\item Evolution of objective function;
\item Local convergence rate of the algorithm;
\end{itemize}

We consider a uniform discretization of the two dimensional bounded domain $\Omega=(0,1) \times (0,1)$ and use the matrices arising from a finite differences discretization of the stationary Bingham model in a pipe. More precisely, we minimize
\begin{equation}
  J(y,u)= \frac{1}{2} |y-1|^2 + \frac{\alpha}{2} |u|^2
\end{equation}
subject to the variational inequality \eqref{eq:vi2},  with $A$ arising from a five point stencil discretization of the Laplacian operator and $\mathbb{K}$ is constructed using centered difference approximations of the first partial derivatives. The mesh size step is set to $h=1/61$. Consequently, the control $u$ is a vector of size $n=61^2$ and the state $y$ is a vector of size $m=61^2$. The Tikhonov parameter $\alpha$ is varied in the range $\alpha \in \{5E-3, 1E-3, 5E-4, 1E-4, 5E-5\}$.

The used parameters for the trust-region algorithm are: $\eta_1=0.25$, $\eta_2=0.75, \gamma_1=0.5, \gamma_2=1.3.$ The initial radius for the algorithm was set to $\Delta_0=10$ and the radius lower bound to $\Delta_{\min}=1E-6$. The second order matrix $H_k$ was built using a standard BFGS approximation. Alternative quasi-Newton updates were not tested, since the BFGS provided satisfactory results. For the fraction of Cauchy decrease condition, we considered $\beta=1$ and $\delta=0.8$. The algorithm starts from the initial constant control $u=10$ and stops whenever $\frac{|u_{k+1}-u_k|}{|u_0|}$ is smaller than a given tolerance, typically set to $1E-4$.

The behaviour of the trust-region algorithm does not depend on the lower-level problem solver. We consider two different type of methods for the Bingham variational inequality. The first one is a semismooth Newton method based on a Huber regularization of the TV term \cite{DelosreyesGonzalez2008}. We tested this algorithm with a regularization parameter $\gamma=1000$. The second algorithm is a primal-dual first order method \cite{TRESKATIS2016115}. In this case no regularization is required, but the number of iterations (and computing time) to reach convergence is much higher. This different behaviour of the lower-level problem solvers, however, does not have an impact on the number of iterations of our TR algorithm. Moreover, both solvers can be combined in order to get an accelerated inexact type algorithm.

\begin{figure}[h]
  \begin{center}
\includegraphics[height=4cm,width=6cm]{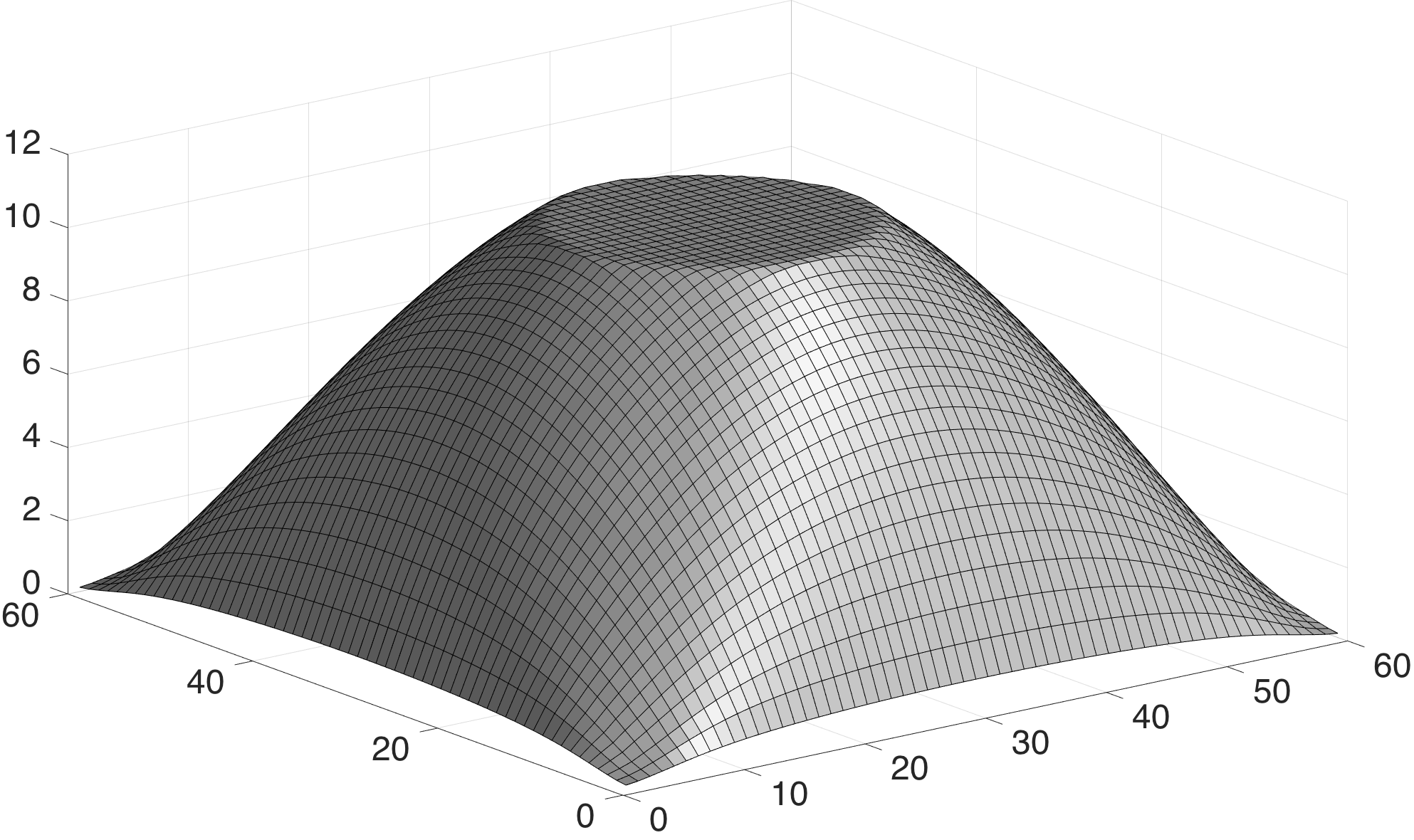} \hfill
\includegraphics[height=4cm,width=6cm]{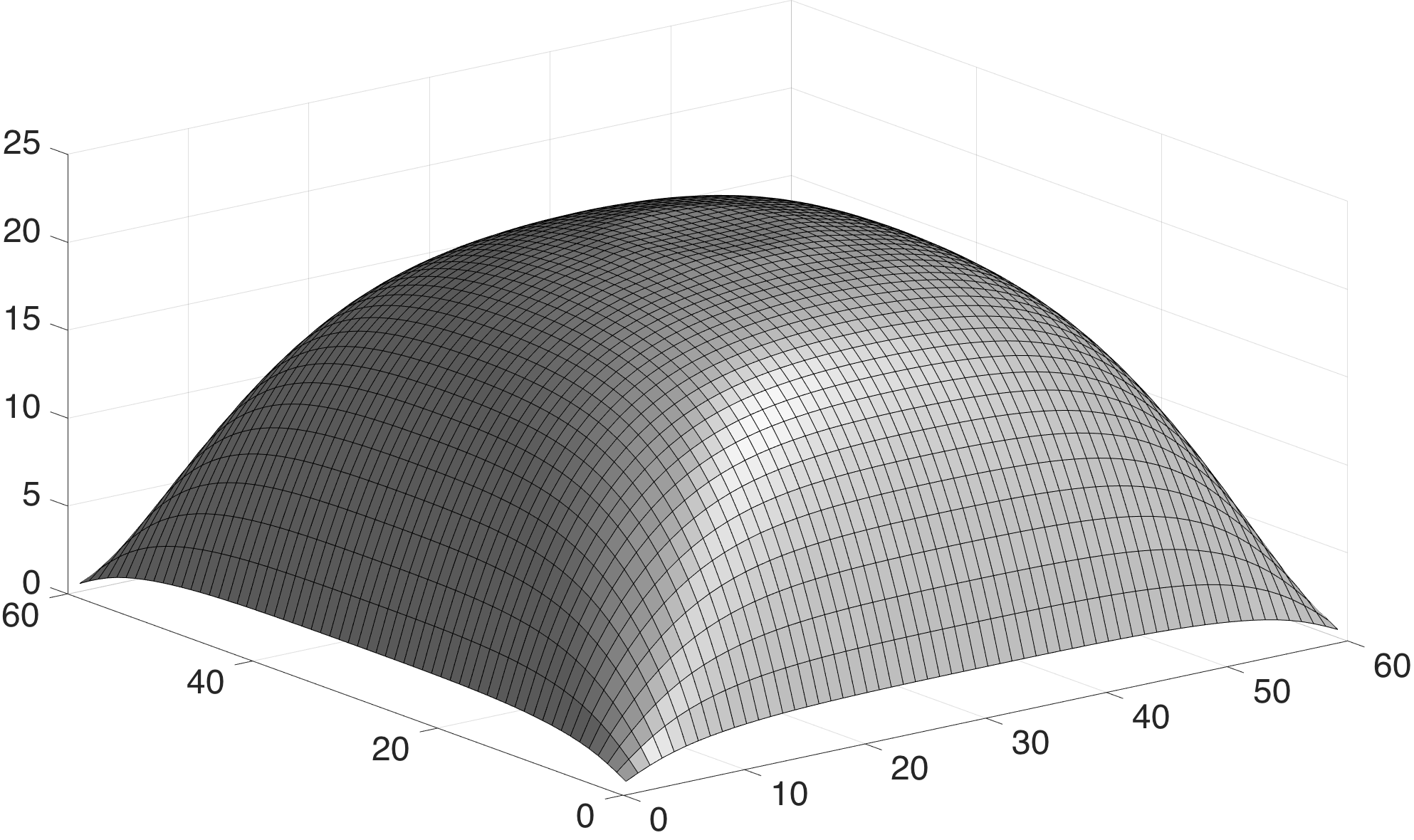}\\
\includegraphics[height=4cm,width=6cm]{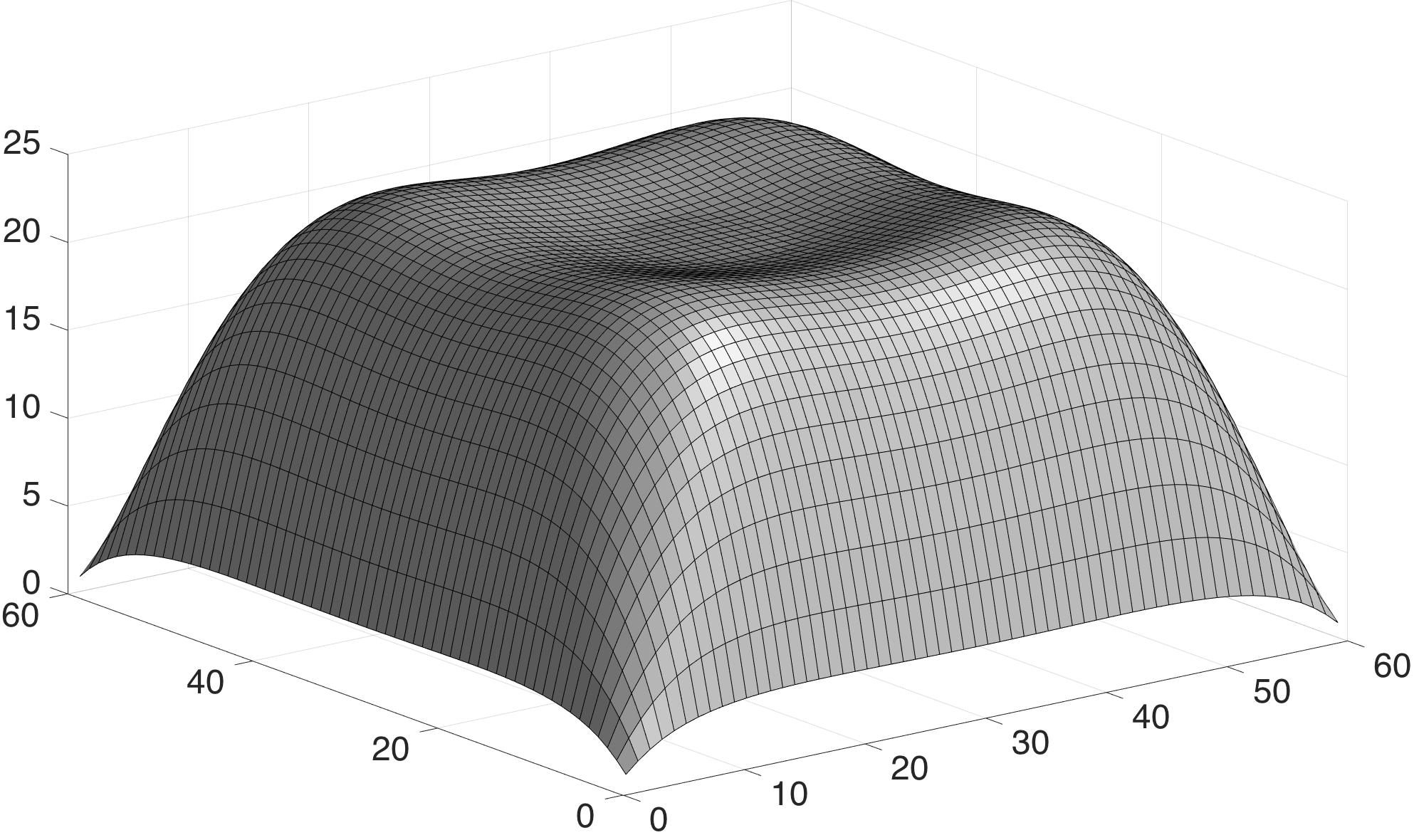} \hfill
\includegraphics[height=4cm,width=6cm]{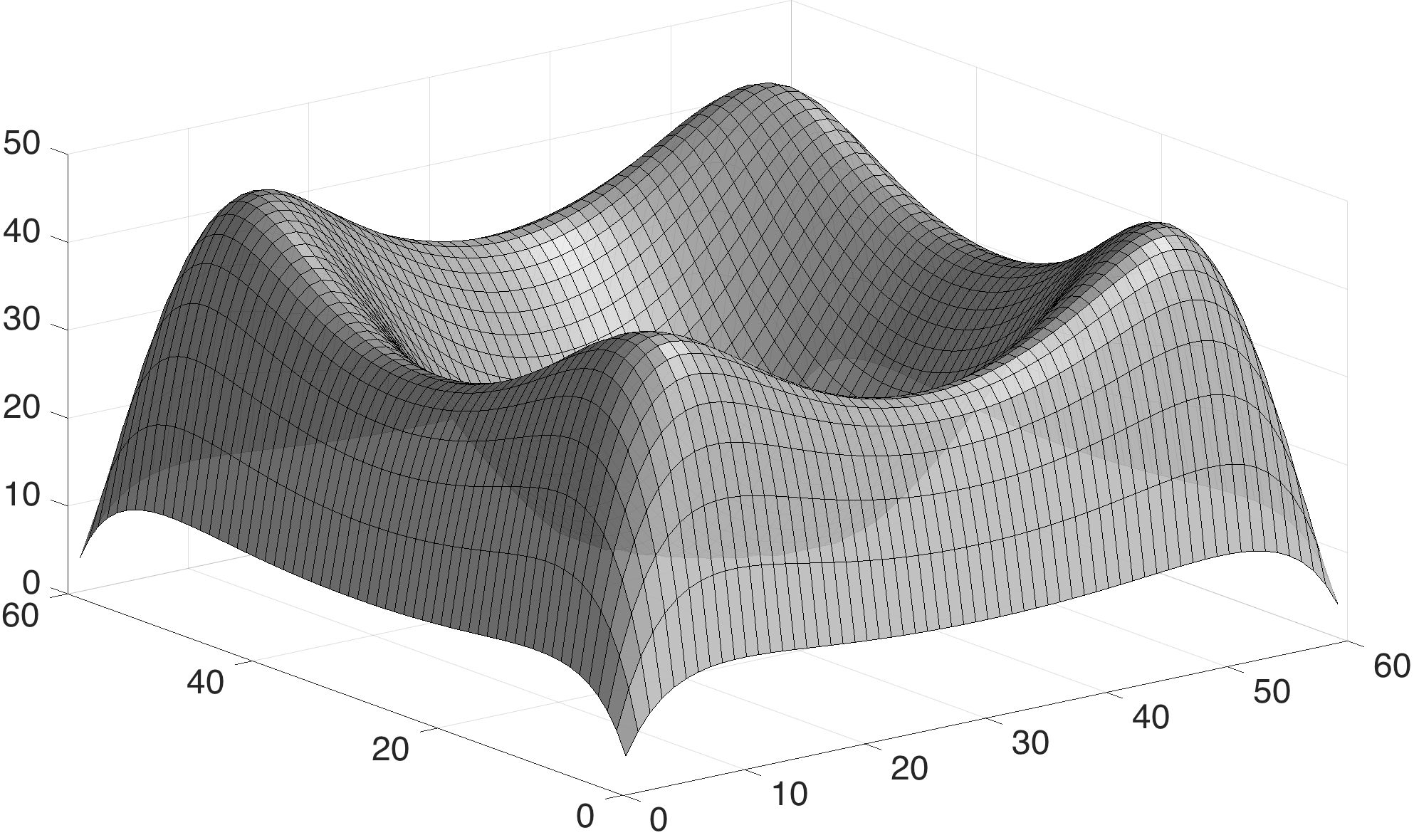}
\caption{Optimal control $u$ for different Tikhonov parameter weights. From the left upper corner to the lower right corner: $\alpha=5E-3$, $\alpha=1E-3$, $\alpha=5E-4$, $\alpha=1E-4$.}
\end{center} \label{figure: controls for different weights}
\end{figure}

\begin{figure}
  \begin{center}
\includegraphics[height=4cm,width=6cm]{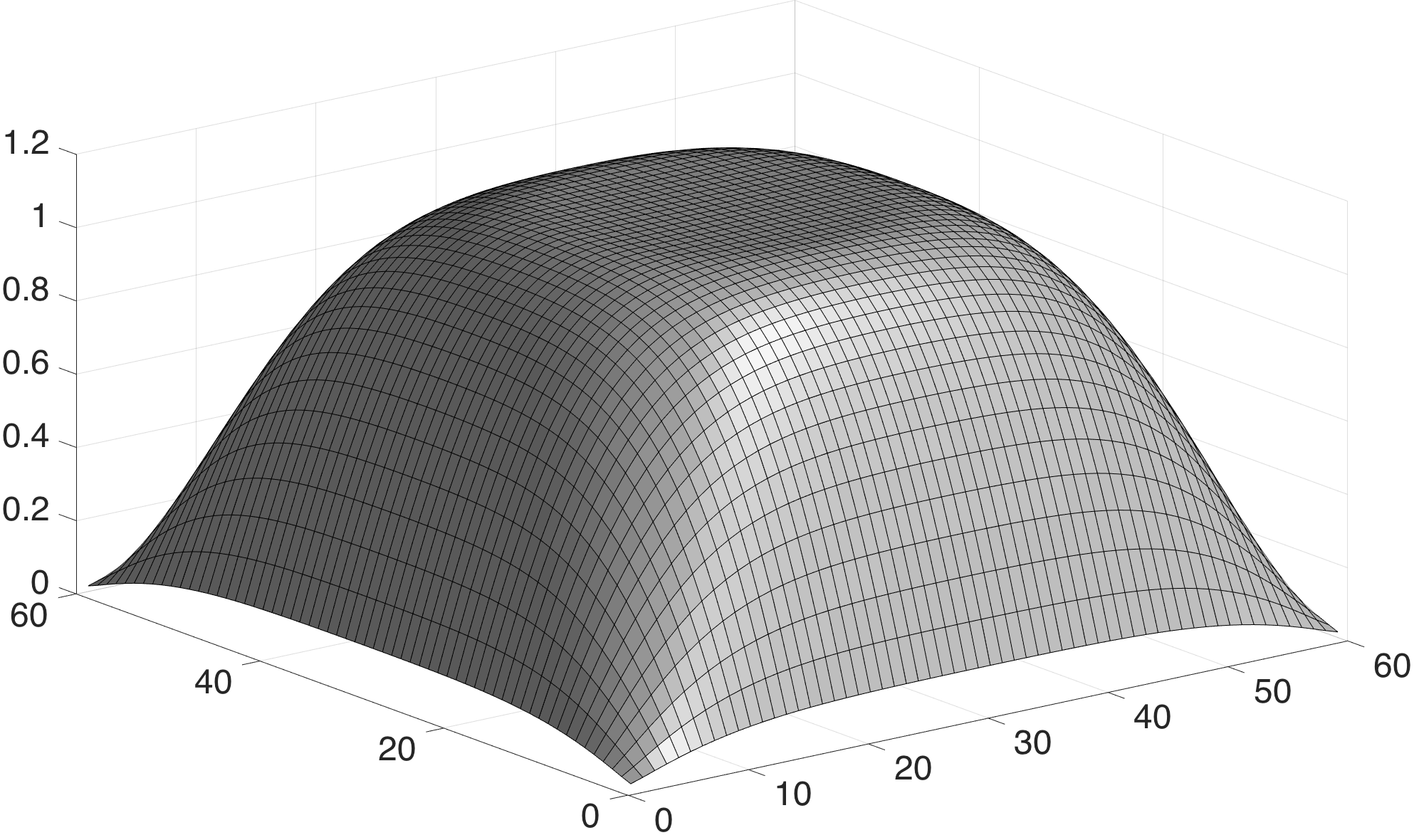} \hfill
\includegraphics[height=4cm,width=6cm]{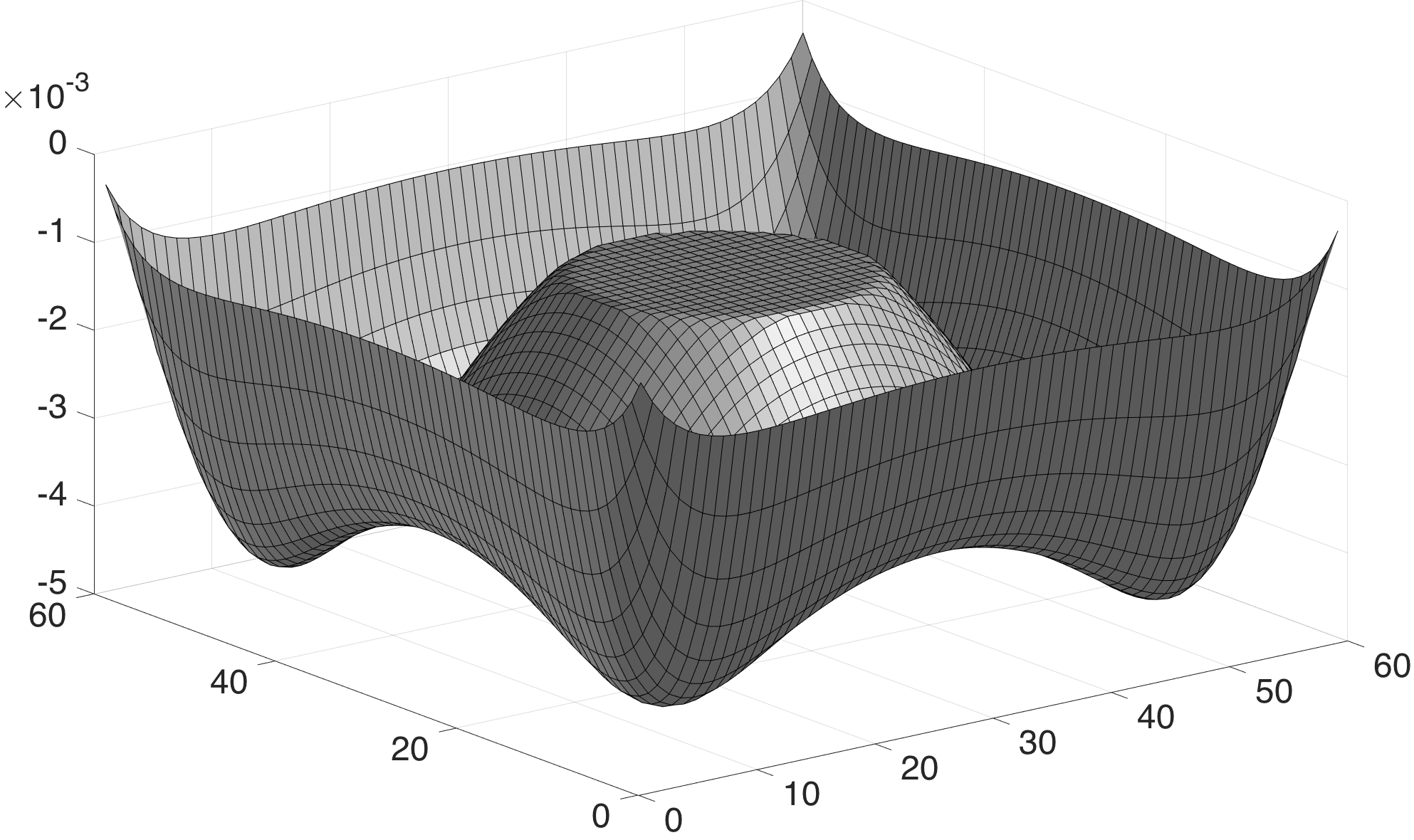}
\caption{Optimal controlled state $y$ and adjoint state $p$ for the control problem with Tikhonov weight $\alpha=1E-4$. Mesh size step $h=1/61$.}
\end{center} \label{figure: optimal state and adjoint}
\end{figure}

Concerning the solution's behaviour, since the desired state is a constant flow velocity equal to one, the optimal control pushes harder close to the boundary as the Tikhonov parameter $\alpha$ becomes smaller. This can be observed from the plots in Figure \ref{figure: controls for different weights}. The computed optimal and adjoint states, for the problem with $\alpha=1E-4$, are depicted in Figure \ref{figure: optimal state and adjoint}, where the resulting nonsmooth structure can be clearly visualized on the adjoint state plot.

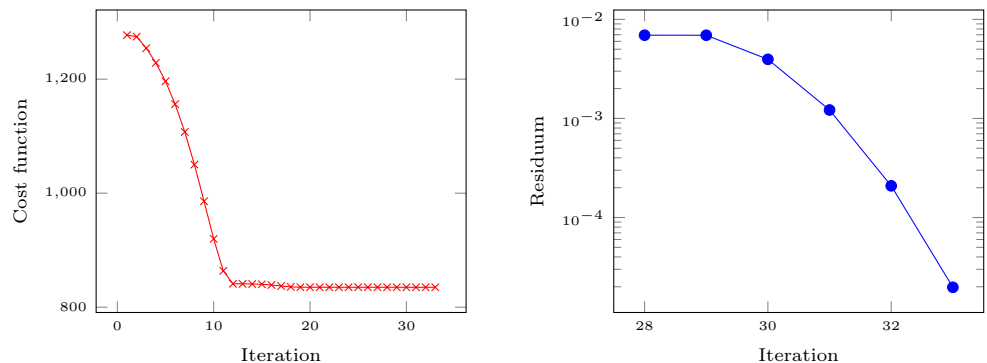
\begin{figure}
\begin{center}
  \begin{tikzpicture}
  	\begin{axis}[
  		xlabel=Iteration,
  		ylabel=Cost function,
      yticklabel style = {font=\tiny},
      ylabel style = {font=\scriptsize},
      xlabel style = {font=\scriptsize},
     xticklabel style = {font=\tiny}]
  	\addplot[color=red,mark=x] coordinates {
    (1,	1277.104634583353)
    (2,	1274.615015207466)
    (3,	1254.312491819300)
    (4,	1228.575949632669)
    (5,	1196.219940174294)
    (6,	1156.133392828301)
    (7,	1107.434898795499)
    (8,	1050.076168483287)
    (9,	985.779452057873)
    (10,	919.760416463701)
    (11,	863.714909338943)
    (12,	840.985180072784)
    (13,	840.975574419701)
    (14,	840.669949495106)
    (15,	840.184668954770)
    (16,	838.959823889841)
    (17,	837.181358459474)
    (18,	835.584564955661)
    (19,	835.020610067957)
    (20,	834.951979144274)
    (21,	834.949569440263)
    (22,	834.949502015696)
    (23,	834.949406569820)
    (24,	834.949106431012)
    (25,	834.948384513970)
    (26,	834.946527010848)
    (27,	834.942228907186)
    (28,	834.933691488585)
    (29,	834.921967187796)
    (30,	834.913624930373)
    (31,	834.911348375435)
    (32,	834.911151639871)
    (33,	834.911146737285)
  	};
  	\end{axis}
  \end{tikzpicture}
\hfill
\begin{tikzpicture}
\begin{semilogyaxis}[
  xlabel=Iteration,
  ylabel=Residuum,
  ylabel style = {font=\scriptsize},
  xlabel style = {font=\scriptsize},
  yticklabel style = {font=\tiny},
 xticklabel style = {font=\tiny}]
\addplot[color=blue,mark=*] coordinates {
  (28,0.006923574607348)
  (29,0.006918740128400)
  (30,0.003952335788742 )
  (31,0.001217345514778)
  (32,0.000208763213798)
  (33,0.000019761178418)
};
\end{semilogyaxis}%
\end{tikzpicture}
\caption{Evolution of the cost funtion along the iterations (left) and residuum in the final 6 iterations of the algorithm (right). Tikhonov parameter $\alpha=5E-4$; mesh size step $h=1/60$.}
\end{center} \label{figure: cost and rate}
\end{figure}

The number of trust-region iterations for different values of the Tikhonov parameter are registered in Table \ref{table: iterations for different Tikhonov}. As expected, as $\alpha$ becomes smaller, the problem is harder to solve and the method requires more iterations. However, the total number of iterations remains small for such difficult problem. Moreover, when looking at the local convergence rate near the solution, a superlinear behaviour can be observed. This is shown in Figure \ref{figure: cost and rate}, together with the evolution of the cost function.
  \begin{table}
  \caption{Number of iterations for different $\alpha$ values}
  \label{table: iterations for different Tikhonov}
  \begin{center}
  \begin{tabular}{|l|c|c|c|c|c|}
    $\alpha$ &5E-3 &1E-3 &5E-4 &1E-4 &5E-5 \\ \hline
    \# iter &24 &29 &33 &55 &58
  \end{tabular}
\end{center}
\end{table}


\section{Conclusions}
The present paper develops a rigorous theoretical framework for analyzing variational inequalities of the second kind involving the discrete total variation. By using a primal-dual reformulation of the VI and a direct quotient analysis, we proved the Bouligand differentiability of the solution operator and provided, for the first time, an explicit and constructive characterization of its Bouligand subdifferential. These theoretical results, aside from being of intrinsic interest, form the cornerstone for deriving sharp optimality conditions, including both Bouligand- and strong-stationarity systems, for discrete optimal control problems governed by total variation-based variational inequalities.  Moreover, the developed framework supports the rigorous design and analysis of trust-region algorithms, which depend critically on a detailed characterization of the solution operator's differentiability properties.

\section*{Data Availability}
No external datasets were used in this study. The code used to implement the trust-region algorithm and reproduce the numerical results is available from the corresponding author upon reasonable request.

\bibliographystyle{plain}
\bibliography{biblio}

\end{document}